%% file: prime.tex
\documentclass[11pt]{amsart}

\usepackage{amssymb}
\usepackage{amscd}
\usepackage{mathrsfs}
\usepackage{graphicx}
\usepackage{color}

\def\mytopsep{3mm}
\newtheoremstyle{myplain}{\mytopsep}{\mytopsep}{\itshape}{0pt}{\bfseries}{.}{3mm}{}
\newtheoremstyle{mydefinition}{\mytopsep}{\mytopsep}{\normalfont}{0pt}{\bfseries}{.}{3mm}{}
\newtheoremstyle{myremark}{\mytopsep}{\mytopsep}{\normalfont}{0pt}{\bfseries}{.}{3mm}{}
\theoremstyle{myplain}
\newtheorem{thm}{Theorem}[section]

\newtheorem{lem}[thm]{Lemma}
\newtheorem{prop}[thm]{Proposition}

\theoremstyle{mydefinition}
\newtheorem{dfn}[thm]{Definition}
\theoremstyle{myremark}
\newtheorem{rem}[thm]{Remark}
\newtheorem{exa}[thm]{Example}
\allowdisplaybreaks[1]

\makeatletter\@addtoreset{equation}{section}\makeatother%
\newcommand\mc{\mathscr}
\newcommand\mb{\mathbf}

\newcommand\mathboxdot{\mathop{\boxdot}\limits}
\newcommand\mathodot{\mathop{\odot}\limits}

\newcommand{\aatop}[2]{\genfrac{}{}{0pt}{}{#1}{#2}}

\newcommand\expcomp[1]{\mc{E} \langle{#1}\rangle}
\DeclareMathOperator{\fix}{fix} \DeclareMathOperator{\aut}{aut}
\DeclareMathOperator{\lcm}{lcm} 
\DeclareMathOperator{\c.t.}{ct}\DeclareMathOperator{\Rec}{Rec}
\DeclareMathOperator{\Par}{Par}\DeclareMathOperator{\Fix}{Fix}

\begin{document}

\title{Prime Graphs and Exponential Composition of Species}
\author{Ji Li}
\address{Department of Mathematics \\
University of Arizona\\
617 N. Santa Rita Ave.\\
P.O.Box 210089\\
Tucson, AZ 85721-0089}
\email{jli@math.arizona.edu}
\date{\today}

\begin{abstract}
In this paper, we enumerate prime graphs with respect to the Cartesian multiplication
of graphs. We use the unique factorization of a connected graph into the product of
prime graphs given by Sabidussi to find explicit formulas for labeled and unlabeled
prime graphs. In the case of species, we construct the exponential composition of species
based on the arithmetic product of species of Maia and M\'endez
and the quotient species, and express the species of connected graphs
as the exponential composition of the species of prime graphs.
\end{abstract}

\maketitle

\section{Introduction}

Under the well-known concept
of Cartesian product of graphs (Definition~\ref{dfn-cartprod}),
prime graphs (Definition~\ref{dfn-primegraph}) are non-trivial
connected graphs that are indecomposable with respect to the Cartesian multiplication.
Hence any connected graph
can be decomposed into a product of prime graphs,
and this decomposition is shown by Sabidussi~\cite{sabidussi} to be unique.

To count labeled prime graphs, we express the Dirichlet exponential
generating series (Definition~\ref{dfn-diriexpseq}) of connected graphs
as the exponential of the Dirichlet exponential generating series of prime graphs.
To count unlabeled prime graphs, we see that
the set of unlabeled connected graphs has a free commutative monoid
structure with its prime set being
the set of unlabeled prime graphs. This free commutative monoid structure
enables us to count unlabeled prime graphs in terms of unlabeled connected graphs.

Eventually, we aim at finding the cycle index of the species of
prime graphs. To be more precise, we want to find the relation
between the species of connected graphs and the species of prime
graphs. To start with, we observe that the species associated to a
graph is isomorphic to the molecular species corresponding to the
automorphism group of this graph. This observation leads to a
relation (Proposition~\ref{graphprod}) between the arithmetic
product of species (Definition~\ref{dfn-arith}), studied by Maia
and M\'endez~\cite{arith}, and the Cartesian product of graphs.
Moreover, a theorem (whose simplified but equivalent version is
given by Proposition~\ref{primepower_aut}) of Sabidussi about the
automorphism groups of connected graphs in terms of the
automorphism groups of their prime factors plays an important
role. We define a new operation, the exponential composition of
species (Definition~\ref{dfn-expcompsp}), which corresponds to the
exponentiation group (Definition~\ref{dfn-wreath}) in the case of
molecular species and is related to the arithmetic product of
species as the composition of species is related to the
multiplication of species. We get a formula
(Theorem~\ref{thm-cycleprime}) expressing the species of connected
graphs as the exponential composition of the species of prime
graphs. The enumeration of the species of prime graphs is
therefore completed by applying the enumeration theorem
(Theorem~\ref{expAF}) for the exponential composition of species,
which is a generalization of an enumeration theorem by Palmer and
Robinson~\cite{palmer} on the cycle index polynomial of the
exponentiation group.

An explicit formula for
the inverse of the exponential composition
would be nice to find, but that problem remains open.

\subsection{Introduction to Species and Group Actions}

The combinatorial theory of species was initiated by Joyal~\cite{joyal, joyal2}.
For detailed definitions
and descriptions about species, readers are referred to~\cite{species-book}.

In short, species are classes of ``labeled structures". More formally,
a \emph{species (of structures)} is a
functor from the category of finite sets with bijections $\mathbb{B}$ to itself.
A species $F$ generates for each finite set $U$ a finite set $F[U]$, which is
 called the set of \emph{$F$-structures} on $U$, and
 for each bijection $\sigma: U \to V$ a bijection $F[\sigma]: F[U]\to
 F[V],$ which is called the \emph{transport of $F$-structures along $\sigma$}.
The symmetric group $\mathfrak{S}_n$ acts on
the set $F[n]=F[\{ 1,2,\dots,n\}]$
by transport of structures. The $\mathfrak{S}_n$-orbits under
this action are called \emph{unlabeled}
$F$-structures of order $n$.

 Each species $F$ is associated with three generating series,
 the \emph{exponential generating series}
 $F(x)=\sum_{n\ge 0} |F[n]|
{x^n}/{n!},$
 the \emph{type generating series}
 $\widetilde{F}(x)=\sum_{n\geq 0} f_n\, x^n,$
 where $f_n$ is the number of unlabeled $F$-structures
 of order $n$, and the \emph{cycle index}
 $$Z_F=Z_F\,(p_1,p_2,\dots)=\sum_{n\geq 0}
 \biggl(\sum_{\lambda \vdash n} \fix \, F[\lambda]\, \frac{p_{\lambda}}
 {z_{\lambda}}
 \biggr),$$
 where $\fix \,F[\lambda]$ denotes the number of
 $F$-structures on $[n]$ fixed by $F[\sigma]$, $\sigma$ is a permutation of $[n]$ with
 cycle type $\lambda$, $p_{\lambda}$ is the power sum symmetric
 function (see Stanley~\cite[p.~297]{ec2}) indexed by the partitions $\lambda$ of
 $n$, and $z_{\lambda}$ is the number of permutations in
 $\mathfrak{S}_n$ that commute with a permutation of cycle type
 $\lambda$.

The following identities (see Bergeron, Labelle, and Leroux~\cite[p.~18]{species-book}) illustrate the
 importance of the cycle index in the theory of species.
 \begin{align*}
  F(x)&=Z_F(x, 0, 0, \dots), \\
  \widetilde{F}(x)&=Z_F(x, x^2, x^3, \dots).
 \end{align*}

For example, let $\mc{G}$ be the species of graphs. Note that by graphs we mean
 simple graphs, that is, graphs without loops or multiple
edges. The cycle index of $\mc{G}$ was given in \cite[p.~76]{species-book}:
\begin{equation}
Z_{\mc{G}}=\sum_{n\ge 0} \biggl( \sum_{\lambda \vdash n} \fix\,
\mc{G}[\lambda] \,\frac{p_{\lambda}}{ z_{\lambda}}\biggr),
\notag
\end{equation}
where
$$ \fix\,\mc{G}[\lambda]=2^{\frac{1}{ 2}\sum_{i,j \ge 1}\gcd(i,\,j)\,c_i(\lambda)c_j(\lambda) - \frac{1}{ 2}
\sum_{k\ge1}(k \bmod 2)\,c_k(\lambda)} ,$$ in which $c_i(\lambda)$
denotes the number of parts of length $i$ in $\lambda$. Let ${\mc{G}^c}$ be the species of connected graphs, and $\mc{E}$  the species of sets.
The observation that every graph is a set of connected graphs
 gives rise to the following species identity:
\begin{equation}
\mc{G}=\mc{E}({\mc{G}^c}), \notag
\end{equation}
which can be read as ``a graph is a set of connected graphs", and
gives rise to the identities
\begin{align}
    {\mc{G}^c}(x) &=\log (\mc{G}(x)), \notag \\
   \widetilde{{\mc{G}^c}}(x)&=\sum_{k\ge 1}\,\frac{\mu(k)}{k}\,\log
  (\widetilde{{\mc{G}}}(x^k)), \notag \\
   Z_{{\mc{G}^c}}&=\sum_{k\ge 1}\,\frac{\mu(k)}{k}\,\log(Z_{\mc{G}}\circ
p_k) \label{cycleindex_conngraph},
\end{align}
where the operator $\circ$ on the right-hand side of~\eqref{cycleindex_conngraph}
denotes the operation of plethysm on symmetric
functions (see Stanley~\cite[p.~447]{ec2}).
For example, we can compute the first several terms of the cycle index of the species
of connected graphs ${\mc{G}^c}$ using Maple:
\begin{align}
    Z_{{\mc{G}^c}} &=  p_1 +  \biggl( \frac{1}{2}\,p_1^2+\,\frac{1}{2}\,p_2 \biggr)
  + \biggl( \frac{1}{3}\,p_3+\,\frac{2}{3}\,p_1^3+p_1p_2 \biggr)
   \notag\\
  & \phantom{=} \ \, +
  \biggl(  \frac{19}{12}\,p_1^4+2p_1^2p_2+\,\frac{5}{4}\,p_2^2+\,\frac{2}{3}\,p_1p_3+\,\frac{1}{2}\,p_4 \biggr)
   \notag \\
  & \phantom{=}\ \, + \biggl( \,{19}{3}\,p_1^3p_2+\,\frac{2}{3}\,p_2p_3+
  \,\frac{91}{15}\,p_1^5+5p_1p_2^2+\,\frac{4}{3}\,p_1^2p_3+\,\frac{3}{5}\,p_5+p_1p_4 \biggr)  \notag\\
      & \phantom{=}\ \,  + \biggl( \frac{1669}{45}\,p_1^6 + \,\frac{91}{3}\,p_1^4p_2
    +\,\frac{38}{9}\,p_1^3p_3 + \,\frac{43}{2}\,p_1^2p_2^2 + 2p_1^2p_4 +\,\frac{8}{3}\,p_1p_2p_3
      \biggr.  \notag\\
    &  \phantom{=}\ \,  + \biggl.  \frac{4}{5}\,p_1p_5 +\,\frac{26}{3}\,p_2^3
    +\,\frac{5}{2}\,p_2p_4+\,\frac{25}{18}\,p_3^2+\,\frac{5}{6}\,p_6   \biggr) + \cdots.
   \label{eq-cycle-conn}
\end{align}



For operations of species, readers are referred to \cite[pp.~1--58]{species-book}
for more detailed definitions of the sum $F_1+F_2$, the product
$F_1  F_2=F_1\cdot F_2$, and the composition $F_1(F_2)=F_1\circ F_2$ of arbitrary species
$F_1$ and $F_2$.

The \emph{quotient species} (see~\cite[p.~159]{species-book}) is defined based
on group actions. It appeared in
\cite{2-tree} and \cite{bousquet} as an important tool in
combinatorial enumeration.
Suppose that a group $A$ acts naturally (see~\cite[p.~393]{species-book})
on a species $F$. The quotient species of $F$ by $A$,
denoted $F/A$, is defined to be such that
for each finite set $U$,  $F/A$-structures on $U$ is the set of
$A$-orbits of $F$-structures on $U$, and
for each bijection $\sigma: U \rightarrow V$, the transport of
structures $(F/A)[\sigma]: F[U]/A \rightarrow F[V]/A$
is induced from the bijection $F[\sigma]$ that sends each $A$-orbit
of the set $F[U]$ to an $A$-orbit of the set $F[V]$.


The notion of \emph{molecular species} plays an important role in the analysis
of species. Roughly speaking, a molecular species
is one that is indecomposable under addition. More precisely,
a species $M$ is molecular~\cite{yeh1, yeh2} if there is
only one isomorphism class of $M$-structures,
i.e., if any two arbitrary $M$-structures are isomorphic.

If $M$ is molecular, then $M$ is concentrated on $n$ for some positive
integer $n$, i.e., $M[U]\neq \emptyset$
if and only if $|U|=n$. If this is the case, then there
is a subgroup $A$ of $\mathfrak{S}_n$ such
that $M$ is isomorphic t to the quotient species of $X^n$,
the species of linear orders on an $n$-element set,
by $A$, i.e., $M=X^n/A$. Furthermore, for $A$ and $B$
two subgroups of $\mathfrak{S}_n$ for some $n$, the molecular species
$X^n/A$ is isomorphic to the molecular species $X^n/B$ if and only if
$A$ and $B$ are conjugate subgroups of $\mathfrak{S}_n$.
In other words, for each positive integer $n$,
we get a bijection $\delta_n$ from the set of
conjugate classes of subgroups of the symmetric group of order $n$ to the set
of molecular species concentrated on the cardinality $n$.
A formal construction for the molecular species $X^n/A$ for
a given subgroup $A$ of $\mathfrak{S}_n$ is given by
Bergeron, Labelle, and Leroux~\cite[p.~144]{species-book}.

P\'olya's \emph{cycle
index polynomial}~\cite[pp.~64--65]{polya} of a subgroup $A$ of $\mathfrak{S}_n$
is defined to be
  $$Z(A) = Z(A; p_1, p_2, \dots, p_n)=
  \frac{1}{|A|} \sum_{\sigma \in A} \prod _{k=1}^n p_k^{c_k(\sigma)},$$
where $c_k(\sigma)$ denotes the number of $k$-cycles in the permutation $\sigma$.

An application of Cauchy-Frobenius Theorem~\cite{robinson} (Lemma~\ref{burnsides})
gives that the cycle index polynomial of $A$
is the same as the cycle index of the molecular species $X^n/A$
(see~\cite[p.~117]{macdonald} Example 7.4):
$$Z(A)=Z_{X^n/A}.$$
This formula illustrates that the cycle index series of species
is a generalization of P\'olya's cycle index polynomial.

\begin{dfn}
\label{dfn-speciesassociatedgraph} An example of molecular species
is the \emph{species associated to a graph}. For each graph $G$ we
assign a species $\mc{O}_G$ to it such that for any finite set
$U$, the set $\mc{O}_G[U]$ is the set of graphs isomorphic to $G$
with vertex set $U$. The species $\mc{O}_G$ is the molecular
species corresponding to the automorphism group of $G$ as a
subgroup of the symmetric group on the vertex set of $G$. We write
$Z(G)$ for the cycle index of the species associated to the graph
$G$, which is the same as the cycle index polynomial of the
automorphism group of $G$. In other words,
$$Z(G)=Z_{\mc{O}_G}=Z(\aut(G)). $$
\end{dfn}

The fact that molecular species are indecomposable under addition
leads to a \emph{molecular decomposition} of any species~\cite[p.~141]{species-book}.
That is, every species of structures $F$ is the sum of its molecular subspecies:
$$F=\sum_{\aatop{M \subseteq F}{ M \text{ molecular}}} M.$$

Let $A$ be a subgroup of $\mathfrak{S}_m$, and let $B$ be a subgroup of
$\mathfrak{S}_n$. We can construct new groups based on $A$ and $B$.

\begin{dfn}
\label{dfn-prodgroup}
The \emph{product group} whose elements are of the form $(a,b)$, where $a\in A$
and $b\in B$, and whose group operation is given by
$(a_1, b_1)\cdot (a_2, b_2)=(a_1a_2, b_1b_2),$
where $a_1$ and $a_2$ are elements of $A$, and $b_1$ and $b_2$ are
elements of $B$, has two group representations, denoted by $A \divideontimes B$ and $A \times B$,
where the group $A \divideontimes B$ acts on the set $[m+n]$ by
\begin{equation}
  (a,b)(i)=\left\{
  \begin{array}{ll}
    a(i), & \text{ if } i \in \{1,2,\dots,m\},\\
    b(i-m)+m, & \text{ if } i \in \{ m+1, m+2, \dots, m+n\},
  \end{array}
  \right.
\end{equation}
and the group $A\times B$ acts on the set $[m]\times [n]$  by
$(a,b)(i, j)=(a(i),b(j)),$
for all $i\in [m]$ and $j \in [n]$.

Therefore, we can identify the group $A \divideontimes B$ with a subgroup
of $\mathfrak{S}_{m+n}$, and the group $A \times B$ with a subgroup of
$\mathfrak{S}_{mn}$.
\end{dfn}

\begin{dfn}
\label{dfn-wreath}
The \emph{wreath product} of $A$ and $B$ has group elements of the form $(\alpha, \tau)$,
where  $\alpha$ is a permutation in $A$ and $\tau$ is a function from
$[m]$ to $B$. The composition of two
elements $(\alpha, \tau)$ and $(\beta, \eta)$ of $B \wr A$ is
given by
$$(\alpha, \tau)(\beta, \eta)=(\alpha\beta, (\tau \circ \beta)
    \eta),$$
where $\beta \in A$ is viewed as a function from $[m]$ to $[m]$,
and $(\tau \circ \beta) \eta$ denotes the point-wise multiplication of $\tau
\circ \beta$ and
$\eta$, both functions from $[m]$ to $B$.

We introduce two group representations of
the wreath product of $A$ and $B$, denoted $B \wr A$ and $B^A$,
which were studied in full detail by Palmer and Robinson~\cite{palmer}.

First, the group $B \wr A$ acts on the set $[m]\times[n]$ by letting
$(\alpha, \tau) (i, j) = (\alpha i, \tau(i) j),$
for all $i \in [m]$ and $j\in [n]$. Hence the group $B \wr A$ can
be identified with a subgroup of $\mathfrak{S}_{mn}$.

Second, the group $B^A$ acts on the set of functions from $[m]$ to $[n]$
by letting $(\alpha, \tau)(f)=g$ for $f:[m]\rightarrow [n]$,
where $g:[m]\rightarrow [n]$ is defined by
$$( (\alpha, \tau) f) (i) =g(i)= \tau(i)( f(\alpha ^{-1} i)),$$
for any $i\in [m]$. We observe that the group $B^A$ can be identified
with a subgroup of $\mathfrak{S}_{n^m}$.
\end{dfn}

 Yeh~\cite{yeh1,yeh2} proved the following species identities:
  \begin{align*}
  \frac{X^m}{A}  \frac{X^n }{B} &=\frac{X^{m+n}}{A \divideontimes
  B},\\
  \frac{X^m}{A}\biggl(\frac{X^n }{ B} \biggr)&=\frac{X^{mn}}{B \wr A}.
  \end{align*}
Note that these results agree with P\'olya's
Theorems~\cite{polya} for the cycle index polynomials of $A \divideontimes B$ and $B \wr A$.
In this paper, we will study the molecular species $X^{mn}/(A \times B)$
(Section~\ref{section_arith}) and the molecular
species $X^{n^m}/B^A$ (Section~\ref{exp-comp}).

\section{Labeled and Unlabeled Prime Graphs}

\subsection{Cartesian Product of Graphs}

For any graph $G$, we let $V(G)$ be the vertex set of $G$, $E(G)$
the edge set of $G$, and  $l(G)=|V(G)|$ the number of vertices in
$G$. Two graphs $G$ and $H$ with the same number of vertices are
said to be \emph{isomorphic}, denoted $G \cong H$, if there exists
a bijection from $V(G)$ to $V(H)$ that preserves adjacency. Such a
bijection is called an \emph{isomorphism} from $G$ to
$H$. In the case when $G$ and $H$ are identical, this bijection
 is called an \emph{automorphism} of $G$. The collection
of all automorphisms of $G$, denoted  $\aut(G)$, constitutes a
group called the \emph{automorphism group} of $G$.
We set $L(G)$ to be the number of
graphs isomorphic to $G$ with vertex set $V(G)$. It is easy to see that
$L(G)={l(G)!}/{\left|\aut(G)\right|}.$
We use the notation $\sum_{i=1}^n G_i=G_1+G_2+\cdots+G_n$ to mean
the disjoint union of a set of graphs $\{G_i \}_{i=1,\dots,n}$.

\begin{figure}[ht]
\begin{center}
\scalebox{1.5}{\input{cart_prod_1.pstex_t}}
\end{center}
\caption[\ \ The Cartesian product of
graphs.]{\label{f-cartprod}The Cartesian product of a graph with
vertex set $\{1,2,3, 4\}$ and a graph with vertex set $\{a,b,c\}$
is a graph with vertex set $\{(i,j)\}$, where $i\in \{1,2,3, 4\}$
and $j\in \{a,b,c\}$.}
\end{figure}
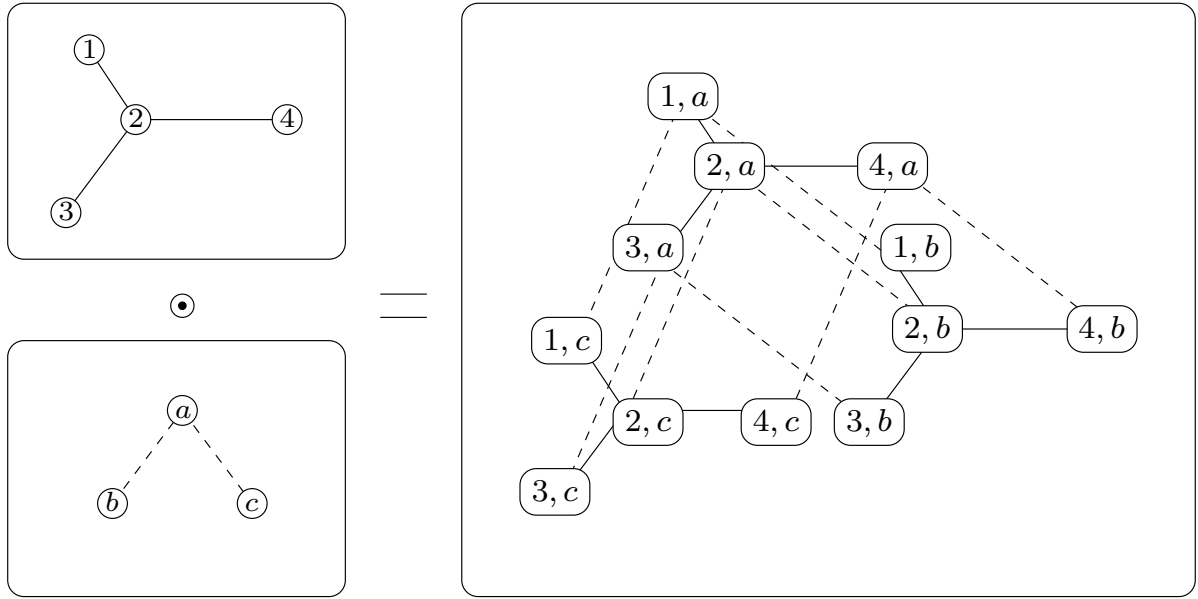

\begin{dfn}\label{dfn-cartprod}
The \emph{Cartesian product} of graphs $G_1$ and $G_2$,
denoted $G_1 \odot G_2$,
as defined by Sabidussi~\cite{sabidussi} under the name the
\emph{weak Cartesian product}, is the graph whose vertex set is
$V(G_1\odot G_2)=V(G_1) \times V(G_2)
=\{ (u,v): u \in V(G_1), v\in V(G_2)\},$
in
which $(u,v)$ is adjacent to $(w,z)$ if
either $u=w$ and $\{ v, z \} \in E(G_2)$ or  $v=z$ and $\{ u, w \} \in E(G_1)$.
\end{dfn}
An example of the Cartesian product of two graphs is given in Figure~\ref{f-cartprod}.

For simplicity and without ambiguity, we call $G_1\odot G_2$ the
\emph{product} of $G_1$ and $G_2$.

It can be verified straightforwardly that the Cartesian multiplication
is commutative and associative up to isomorphism. We denote by $G^n$ the
Cartesian product of $n$ copies of $G$.

\begin{dfn}\label{dfn-primegraph}
A graph $G$ is \emph{prime} with respect to Cartesian
multiplication if $G$ is a connected graph with more than one
vertex such that  $G \cong H_1 \odot H_2$ implies that either $H_1$
or $H_2$ is a singleton vertex.

Two graphs $G$ and $H$ are called \emph{relatively prime} with
respect to Cartesian multiplication, if and only if $G \cong G_1
\odot J$ and $H \cong H_1 \odot J$  imply that $J$ is a singleton
vertex.
\end{dfn}

We denote by ${\mc{P}}$ the species of prime graphs.
We see from Definition~\ref{dfn-primegraph} that any non-trivial connected graph
can be decomposed into a product of prime graphs. Sabidussi~\cite{sabidussi} proved
that such a prime decomposition is unique up to isomorphism.

The automorphism groups of the Cartesian product of a set of graphs
was studied by Sabidussi~\cite{sabidussi} and Palmer~\cite{exponentiation}.
For example, Sabidussi proved that the automorphism group of the
disjoint union of a set of graphs is isomorphic to the automorphism
group of the Cartesian product of these graphs. Sabidussi also showed that
the automorphism group of the Cartesian prduct of the disjoint union of
two relatively prime graphs is the product
of the automorphism group of these two graphs.

\subsection{Labeled Prime Graphs}

In this section all graphs considered are connected.

Sabidussi gave an important formula
about the automorphism group of a connected
graph using its prime factorization:
If $G$ is a connected graph with prime factorization
$$G \cong P_1^{s_1} \odot P_2^{s_2} \odot \cdots \odot
P_k^{s_k},$$
where for $r=1,2,\dots, k$, all $P_r$ are distinct
prime graphs, and all $s_r$ are positive integers, then
$$\aut(G)
\cong \prod_{r=1}^k \aut(P_r^{s_r})
 \cong \prod_{r=1}^k
 \aut(P_r)^{\mathfrak{S}_{\, s_r}}.$$
  Note that the $P_r^{s_r}$, for $r=1,2,\dots,k$, are pairwise relatively prime.
Since the automorphism group of the Cartesian product of the disjoint union of
two relatively prime graphs is the product
of the automorphism groups of the graphs, we see that Sabidussi's formula reduces
equivalently to the following proposition:

\begin{prop}
\label{primepower_aut}
\emph{(Sabidussi~\cite{sabidussi})}
  Let $P$ be a prime graph, and let $k$ be a nonnegative integer.
  Then the automorphism group of $P^k$ is the exponentiation group
  $\aut(P)^{\mathfrak{S}_k}$, i.e.,
  $$\aut(P^k)=\aut(P)^{\mathfrak{S}_k}.$$
In particular,  $$|\aut(P^k)|=|\aut(P)^{\mathfrak{S}_k}|
= k!\cdot |\aut(P)|^k.$$
\end{prop}

\begin{dfn}\label{dfn-diriexpseq}
The \emph{Dirichlet exponential generating series} for a
sequence of numbers $\{a_n\}_{n \in \mb{N}}$ is defined by
$\mathop{\sum}\limits_{n\ge 1} {a_n}/{(n! \, n^s)}.$

Multiplication of Dirichlet exponential generating series is given by
\begin{equation}
\biggl( \sum_{n\ge 1} \frac{a_n}{n!\,n^s} \biggr) \biggl(  \sum_{n
\ge 1} \frac{b_n}{n!\,n^s}\biggr) =\sum_{n\ge 1}
\frac{c_n}{n!\, n^s}, \notag
\end{equation}
where
\begin{equation}
c_n=
\sum_{k|n} \biggl\{ \aatop{n}{k} \biggr\} \,a_k b_{n/k} =
\sum_{k|n}  \frac{n!}{k!\, (n/k)!}  \,a_k b_{n/k}.
\notag
\end{equation}

The \emph{Dirichlet exponential generating function for a species
$F$} with the restriction $F[\emptyset]=\emptyset$ is defined by
$$\mathfrak{D}(F)=\sum_{n\ge 1} \frac{|F[n]|}{n! \, n^s}.$$
\end{dfn}

The \emph{Dirichlet exponential generating function for a graph
$G$} is defined by
$$\mathfrak{D}(G)=\frac{L(G)}{l(G)!\cdot l(G)^s},$$
where $L(G)$ is the number of graphs isomorphic to $G$ with vertex
set $V(G)$, and $l(G)$ is the number of vertices of $G$. In other
words, $$\mathfrak{D}(G)=\mathfrak{D}(\mc{O}_G),$$ where
$\mc{O}_G$ is the species associated to a graph defined by
Definition~\ref{dfn-speciesassociatedgraph}. Recall that
$$L(G)=\frac{l(G)!}{|\aut(G)|}.$$ Therefore,
$$\mathfrak{D}(G)=\frac{1}{|\aut(G)|\cdot l(G)^s}.$$

\begin{exa}
Let ${\mc{P}}$ be the species of prime graphs, let ${\mc{G}^c}$ be the
species of connected graphs, let $\mathbb{C}$ be the set of unlabeled connected graphs,
and let $\mathbb{P}$ be the set of unlabeled prime graphs. Then $\mathfrak{D}({\mc{G}^c})$ and
$\mathfrak{D}({\mc{P}})$ are the Dirichlet exponential generating functions
for these two species, respectively:
\begin{align*}
 \mathfrak{D}({\mc{G}^c})=
 \sum_{n\ge1}\frac{|{\mc{G}^c}[n]|}{n!\,n^s}
 =\sum_{G \in \mathbb{C}} \mathfrak{D}(G),\qquad
  \mathfrak{D}({\mc{P}}) =
 \sum_{n\ge1}\frac{|{\mc{P}}[n]|}{n!\,n^s}
 =\sum_{P \in \mathbb{P}} \mathfrak{D}(P).
\end{align*}
\end{exa}

Propositions~\ref{dirichexpsp} and~\ref{graphprod} lead
straightforwardly to the following lemma.

\begin{lem}\label{prodrelativeprime0}
  Let $G_1$ and $G_2$ be relatively prime graphs. Then
  \begin{equation}
    \mathfrak{D}({G_1 \odot G_2}) = \mathfrak{D}({G_1})\,
    \mathfrak{D}({G_2})
    \label{prodrelativeprime1}
  \end{equation}
\end{lem}

\begin{lem}\label{primepows0}
  Let $P$ be any prime graph. Let $T$  be the set of all nonnegative integer powers of $P$, i.e.,
  $T=\mathop{\cup}\limits_{k \ge 0}P^k$. Then
  the Dirichlet exponential generating functions for $T$ and $P$
  are related by
  \begin{equation}
    \mathfrak{D}(T)=\exp (\mathfrak{D}(P)).
    \label{primepows}
  \end{equation}
\end{lem}

\begin{proof}
We start with
$$\mathfrak{D}(P)=\frac{L(P)}{l(P)!\cdot l(P)^s}=\frac{1}{|\aut(P)|\cdot l(P)^s}.$$
It follows from Proposition~\ref{primepower_aut} that
$$L(P^k) =\frac{l(P^k)!}{|\aut(P^k)|}=\frac{l(P^k)!}{k! \cdot\left| \aut(P) \right|^k },$$
and that
$$
\mathfrak{D}({P^k})=\frac{L(P^k)}{l(P^k)! \cdot l(P^k)^s}
=\frac{1}{k! \cdot\left| \aut(P) \right|^k \cdot l(P)^{ks}}
=\frac{\mathfrak{D}({P})^k}{k!}.
$$
Summing up on $k$, we get
$$\mathfrak{D}(T)=\sum_{k\ge0}\frac{\mathfrak{D}({P})^k}{k!}
= \exp (\mathfrak{D}(P)).$$
\end{proof}

\begin{thm}
\label{thm-degfconnprime} For $\mathfrak{D}({\mc{G}^c})$ and
$\mathfrak{D}({\mc{P}})$, we have
$$\mathfrak{D}({\mc{G}^c})=\exp\, (\mathfrak{D}({\mc{P}})).$$
\end{thm}

\begin{proof}

 Lemma~\ref{prodrelativeprime0} gives that the Dirichlet exponential
generating function of a product of two relatively prime graphs is
the product of the Dirichlet exponential generating functions of
the two graphs. Since the operation of Cartesian
product on graphs is associative up to isomorphism, it
follows that if we have a set of pairwise relatively prime graphs
$\{G_i\}_{i=1, 2, \dots, r}$, and let $G =\mathodot_{i=1}^r
G_i$, then
\begin{equation}
\mathfrak{D}({G})=\prod_{i=1}^r \mathfrak{D}({G_i}).
\label{dirichsetrelprime}
\end{equation}

Now according to the definition of the Dirichlet exponential
generating function for graphs, we get
\begin{align*}
\mathfrak{D}({{\mc{G}^c}}) &= \sum_{G \in \mathbb{C}} \mathfrak{D}(G)
=\prod_{P \in \mathbb{P}}
\mathfrak{D}\biggl({\sum_{k\ge0}P^k}\biggr)
=\prod_{P \in \mathbb{P}}\exp(\mathfrak{D}(P))\\
&=\exp \Big(\sum_{P\in\mathbb{P}} \mathfrak{D}(P) \Big)
=\exp(\mathfrak{D}({\mc{P}})).
\end{align*}
\end{proof}

It is well-known that the exponential generating series of the species
of connected graphs $\mc{G}^c$ is
\begin{align*}
{\mc{G}^c}\,(x) & =  \sum_{n \ge 1} |\mc{G}^c[n]|\,\frac{x^n}{n!}
  = \log \biggl( \sum_{n\ge 1} 2^{\;(\aatop{n}{2})} \,\frac{x^n}{n!}  \biggr) \\
 &=\frac{x }{1!}\,+\, \frac{x^2}{2!}\,+4\,\frac{x^3}{ 3!}\,+{38}\,\frac{x^4
  }{ 4!}\,+ 728\,\frac{x^5}{ 5!}\,+{ 26704}\,\frac{x^6}{ 6!}\,+
  {1866256}\,\frac{x^7}{7!}\, \\
  & \phantom{=}\ \,
   + {251548592}\,\frac{x^8}{ 8!}\,+{66296291072}\,\frac{x^9}{ 9!}\,+\dots.
\end{align*}

We obtain $\mathfrak{D}({\mc{G}^c})$ by replacing  $x^n$ with $n^{-s}$ for
each $n$ in the above expression:
\begin{align*}
\mathfrak{D}({\mc{G}^c}) &=  \sum_{n\ge1}|{\mc{G}^c}[n]|\,\frac{1}{n!\,n^s}\\
&= \frac{1}{1!\,1^s}\,+\, \frac{1}{2!\,2^s}\,+4\, \frac{1}{3!\,3^s}\,+{38}\,\frac{1}{4!\,4^s}\,
+ 728\,\frac{1}{5!\,5^s}\,+{ 26704}\,\frac{1}{6!\,6^s}\,+
  {1866256}\,\frac{1}{7!\,7^s}\, \\
& \phantom{=}\ \, + {251548592}\,\frac{1}{8!\,8^s}\,
  +{66296291072}\,\frac{1}{9!\,9^s}\,+\dots.
\end{align*}

 Theorem~\ref{thm-degfconnprime}
gives a way of counting labeled prime graphs by writing
$$\mathfrak{D}({\mc{P}})=\log \mathfrak{D}({\mc{G}^c}).$$

For example, we write down the first terms of $\mathfrak{D}({\mc{P}})$ as follows:
\begin{align*}
  \mathfrak{D}({\mc{P}}) &=
  \frac{1}{2!\,2^s}\,+4\, \frac{1}{3!\,3^s}\,+{35}\,\frac{1}{4!\,4^s}\,
+ 728\,\frac{1}{5!\,5^s}\,+{ 26464}\,\frac{1}{6!\,6^s}\,+
  {1866256}\,\frac{1}{7!\,7^s}\, \\
& \phantom{=}\ \, + {251518352}\,\frac{1}{8!\,8^s}\,
  +{66296210432}\,\frac{1}{9!\,9^s}\,+\dots.
\end{align*}

\subsection{Unlabeled Prime Graphs}\label{unlabelconn}

In this section all graphs considered are unlabeled and connected.

\begin{dfn}
\label{dfn-Dirichlet}
  The \emph{(formal) Dirichlet series} of a sequence $\{{a_n}\}_{n=1,2,\dots,\infty}$ is defined
  to be
$\sum_{n=1}^\infty {a_n /n^s}$.

The multiplication of Dirichlet series
 is given by
$$\sum_{n\ge 1} \frac{a_n }{ n^s} \cdot \sum_{m \ge 1} \frac{b_n }{n^s}
=\sum_{n\ge 1}\biggl( \sum_{k|n}a_k b_{n/k} \biggr) \frac{1}{n^s}.$$
\end{dfn}

\begin{dfn}
 A \emph{monoid} is a semigroup with a unit. A \emph{free
 commutative monoid} is a commutative monoid $M$ with a set of primes $P \subseteq M$ such
that each element $m\in M$ can be uniquely decomposed into a
product of  elements in $P$ up to rearrangement. Let $M$ be a free commutative monoid. We get a \emph{monoid
algebra} $\mb{C}M$, in which the elements are all formal sums $
\sum_{m\in M} c_m m,  \text{ where } c_m \in \mb{C},$ with
addition and multiplication defined naturally.
For each $m\in M$, we associate a \emph{length} $l(m)$
that is compatible with the multiplication in $M$. That is, for
any $m_1, m_2\in M$, we have $l(m_1)l(m_2)=l(m_1m_2)$.
\end{dfn}

Let $M$ be a free commutative monoid with prime set $P$. The
following identity holds in the monoid algebra $\mb{C}M$:
\begin{equation}\notag
  \sum_{m\in M} m =\prod_{p\in P} \frac{1}{1-p}.
\end{equation}
Furthermore, we can define a homomorphism from $M$ to the ring of
Dirichlet series under which each $m\in M$ is sent to ${1 /
l(m)^s}$, where $l$ is a length function of $M$. Therefore,
\begin{equation}
  \sum_{m\in M} \frac{1}{l(m)^s}=\prod_{p\in P}
  \frac{1}{1-l(p)^{-s}}. \notag
\end{equation}

Recall that $\mathbb{C}$  is the set of unlabeled connected graphs
under the operation of Cartesian product. The unique factorization
theorem of Sabidussi gives $\mathbb{C}$ the structure of a
commutative free monoid with a set of primes $\mathbb{P}$, where
$\mathbb{P}$ is the set of unlabeled prime graphs. This is saying
that every element of $\mathbb{C}$ has a unique factorization of
the form $b_1^{e_1}b_2^{e_2} \cdots b_k^{e_k}$, where the $b_i$
are distinct primes in $\mathbb{P}$.
Let $l(G)$, the number of vertices in $G$, be a length function
for $\mathbb{C}$. We have the following proposition.

\begin{prop}
For $\mathbb{C}$ and $\mathbb{P}$, we have
\begin{equation}\label{monoid}
\sum_{G\in \mathbb{C}} \frac{1}{l(G)^s}=\prod_{P \in
\mathbb{P}}\frac{1}{1-l(P)^{-s}}. \notag
\end{equation}
\end{prop}

The enumeration of prime graphs was studied by Rapha\"{e}l
Bellec~\cite{raphael}. We use Dirichlet series to
count unlabeled connected prime graphs.

\begin{thm}\label{cndnbn}
Let $\widetilde{c}_n$ be the number of unlabeled connected graphs on $n$
vertices, and let $b_m$ be the number of unlabeled prime graphs on
$m$ vertices. Then we have
\begin{equation}
\sum_{n\ge 1} \frac{\widetilde{c}_n}{ n^s}=\prod_{m\ge 2}
\frac{1}{(1-m^{-s})^{b_m}}. \label{unlabelprime}
\end{equation}

Furthermore, if we define numbers $d_n$ for positive integers $n$
by
\begin{equation}\label{dncn}
\sum_{n\ge 1} \frac{d_n}{ n^s} = \log \sum_{n\ge 1} \frac{\widetilde{c}_n}{ n^s},
\end{equation}
then
\begin{equation}\label{dncn2}
d_n=\sum_{m^l=n} \frac{b_m}{ l},
\end{equation}
 where the sum is
over all pairs $(m,l)$ of positive integers with $m^l=n$.
\end{thm}

The proof of Therem~\ref{cndnbn} follows Remark~\ref{rem_cndnbn} and Proposition~\ref{prop_cndnbn} below.

\begin{rem}\label{rem_cndnbn}
In what follows, we introduce an interesting recursive formula for
computing $d_n$. To start with, we differentiate both sides of
Equation~\eqref{dncn} with respect to $s$ and simplify. We get
that
$$
\sum_{n\ge 2} \log n  \frac{\widetilde{c}_n }{ n^s} = \biggl(\sum_{n\ge 1}
\frac{\widetilde{c}_n}{ n^s} \biggr)
 \biggl( \sum_{n\ge 2} \log n  \frac{d_n }{ n^s}\biggr),
$$
 which gives
\begin{equation}\label{dncn1}
\widetilde{c}_n \log n =\sum_{ml=n}  \widetilde{c}_m  d_l \log l .
\end{equation}
Since $\widetilde{c}_1$ is the number of connected graphs on $1$ vertex,
$\widetilde{c}_1=1$. It follows easily from Equation~\eqref{dncn1} that
$d_p=\widetilde{c}_p$ when $p$ is a prime number. Therefore, if $p$ is a prime
number, $b_p=d_p=c_p$. This fact can be seen directly, since a
connected graph with a prime number of vertices is a prime graph.

Rapha\"{e}l Bellec
used Equation~\eqref{dncn1} to find formulae
for $d_n$ where $n$ is a product of two different primes or a
product of three different primes:

If $n=pq$ where $p\neq q$,
\begin{equation}\label{2.4.7}
d_n=\widetilde{c}_n -\widetilde{c}_p\widetilde{c}_q;
\end{equation}

If $n=pqr$ where $p, q$ and $r$ are distinct primes,
\begin{equation}\label{2.4.8}
d_n=\widetilde{c}_n +2\widetilde{c}_p\widetilde{c}_q\widetilde{c}_r-\widetilde{c}_p\widetilde{c}_{qr}
-\widetilde{c}_q\widetilde{c}_{pr}-\widetilde{c}_r\widetilde{c}_{pq}.
\end{equation}

In fact, Equations~\eqref{2.4.7} and~\eqref{2.4.8} are special
cases of the following proposition.

\begin{prop}\label{prop_cndnbn}
Let $d_n, \widetilde{c}_n$ be defined as above. Then we have
\begin{equation}
\label{dncn3}
  d_n  = \widetilde{c}_n - \frac{1}{ 2}\sum_{n_1n_2=n} \widetilde{c}_{n_1}
  \widetilde{c}_{n_2} +  \frac{1}{3}\sum_{n_1n_2n_3=n}\widetilde{c}_{n_1}
  \widetilde{c}_{n_2}\widetilde{c}_{n_3}-\dots.
\end{equation}
\end{prop}

\begin{proof}
We can use the identity
$$\log(1+x)=x-\,\frac{1}{ 2}\,x^2+\,\frac{1}{ 3}\,x^3-\,\frac{1}{ 4}x^4+\dots$$
to compute from Equation~\eqref{dncn} that
\begin{align*}
\sum_{n\ge 1} \frac{d_n}{ n^s}
=& \log \biggl( 1+ \sum_{n\ge 2} \frac{\widetilde{c}_n }{ n^s}\biggr)\\
=& \sum_{n\ge 2} \frac{\widetilde{c}_n }{n^s} - \,\frac{1}{2} \, \biggl(
\sum_{n\ge 2}\frac{\widetilde{c}_n }{ n^s}\biggr)^2+\,\frac{1}{3}\, \biggl(
\sum_{n\ge 2}\frac {\widetilde{c}_n }{ n^s}\biggr)^3-\dots.
\end{align*}
Equating coefficients of $n^{-s}$ on both sides, we get equation
\eqref{dncn3} as a result.
\end{proof}
\end{rem}

\begin{proof}[Proof of Theorem~\ref{cndnbn}]
We start with
\begin{equation}
\sum_{m} \frac{1}{l(m)^s}=\prod_{p}\frac{1}{1-l(p)^{-s}},
\end{equation}
where the left-hand side is multiplied over all  connected graphs, and
the right-hand side is summed over all  prime graphs. Regrouping
the summands on the left-hand side with respect to the number of
vertices in $m$, we get the left-hand side of Equation~\eqref{unlabelprime}.
Regrouping the factors on the right-hand
side with respect to the number of vertices in $p$, we get the
right-hand side of Equation~\eqref{unlabelprime}.


Taking the logarithm of both sides of Equation~\eqref{unlabelprime}, we get
\begin{align*}
\log \sum_{n\ge 1} \frac{\widetilde{c}_n}{ n^s} &= \log \prod_{m\ge 2}
\frac{1}{(1-m^{-s})^{b_m}} = \sum_{m\ge 2} b_m \log \frac{1}{1-m^{-s}}\\
&= \sum_{m\ge2} \biggl ( b_m  \sum_{l \ge 1}
\frac{m^{-sl}}{l} \biggr) = \sum_{m \ge 2,\ l\ge 1}\frac{b_m}{l\ m^{sl}},
\end{align*}
and Equation~\eqref{dncn2} follows immediately.
\end{proof}

Next, we will  compute the numbers $b_n$ in terms of the numbers
$d_n$ using the following lemma.

\begin{lem}\label{dnpn}
Let $\{D_i\}_{i=1,\dots}$ and $\{J_i\}_{i=1,\dots}$ be sequences
of numbers satisfying
\begin{equation}
D_k=\sum_{l|k}\frac {J_{k/l}}{l}, \label{DkPk}
\end{equation}
and let $\mu$ be the M\"{o}bius function. Then we have
\begin{equation*}
J_k=\frac{1}{ k}\sum_{l|k} \mu\biggl(\frac{k}{ l}\biggr)\, l D_l.
\end{equation*}
\end{lem}

\begin{proof}
Multiplying by $k$ on both sides of Equation~\eqref{DkPk} , we get
$$k D_k =\sum_{l|k}\frac{k}{l}\,  J_{k/l} =\sum_{l|k} l J_l.$$
Applying the M\"{o}bius inversion formula, we get
$$kJ_k=\sum_{l|k} \mu\biggl(\frac{k}{l}\biggr) l D_l.$$
Therefore,
$$J_k=\frac{1}{ k}\sum_{l|k} \mu \biggl(\frac{k}{l}\biggr)\, l D_l.$$
\end{proof}

Given any natural number $n$, let $e$ be the largest number such
that $n=r^e$ for some $r$. Note that $r$ is not a power of a
smaller integer. We let $D_k=d_{r^k}, J_k=b_{r^k}$. It follows
that Equation~\eqref{dncn2} is equivalent to Equation~\eqref{DkPk}.

\begin{thm}
\label{thm-unlabelprime}
  For any natural number $n$, let $e, r$ be as described in above.
  Then we have
  \begin{equation*}
    b_n=\frac{1}{ e}\, \sum_{l | e} \mu \biggl(\frac{e}{ l} \biggr) l
    d_{r^e}.
  \end{equation*}
\end{thm}

\begin{proof}
  The result follows straightforwardly from Lemma~\ref{dnpn}.
\end{proof}

Table~\ref{t-lpg} in the Appendix gives the numbers of labeled and unlabeled prime graphs
with no more than $16$ vertices.

\section{Exponential Composition of Species}

\subsection{Arithmetic Product of Species}
\label{section_arith}
The arithmetic product was studied by Maia and M\'{e}ndez~\cite{arith}.
The \emph{arithmetic product} of two molecular species
$X^m/A$ and $X^n/B$, where $A$ is a subgroup of $\mathfrak{S}_m$ and $B$
is a subgroup of $\mathfrak{S}_n$, can be defined to be the molecular species
$X^{mn}/(A \times B)$, where $A \times B$ is the group representation of the product
group of $A$ and $B$ acting on the set $[m]\times [n]$ (Definition~\ref{dfn-prodgroup}).

In order to define the arithmetic
product of general species, Maia and M\'{e}ndez developed a
decomposition of a set, called a \emph{rectangle}.

\begin{dfn}\label{dfn-rectangle}
  Let $U$ be a finite set. A \emph{rectangle}
  on $U$ of \emph{height} $a$ is a pair $(\pi_1,\pi_2)$ such that $\pi_1$
  is a partition of $U$ with
  $a$ blocks, each of size $b$, where $|U|=ab$, and $\pi_2$ is a partition of $U$ with $b$ blocks, each of
  size $a$, and if $B$ is a block of $\pi_1$ and $B'$ is a block of
  $\pi_2$ then $|B \cap B'|=1$.

  A  \emph{$k$-rectangle} on $U$ is a $k$-tuple of
  partitions $(\pi_1, \pi_2, \dots, \pi_k)$ such that

  i) for each $i \in [k]$, $\pi_i$ has $a_i$ blocks, each of
    size $|U|/a_i$, where $|U|=\prod_{i=1}^k a_i$.

  ii) for any $k$-tuple $(B_1, B_2, \dots, B_k)$, where
    $B_i$ is a block of $\pi_i$ for each $i \in [k]$, we have $|
    B_1 \cap B_2 \cap \cdots \cap B_k|=1$. See Figure~\ref{f-rec3_detail_alt}
    for a $3$-rectangle $(\pi_1, \pi_2, \pi_3)$ represented by a $3$-partite
    graph.
\end{dfn}
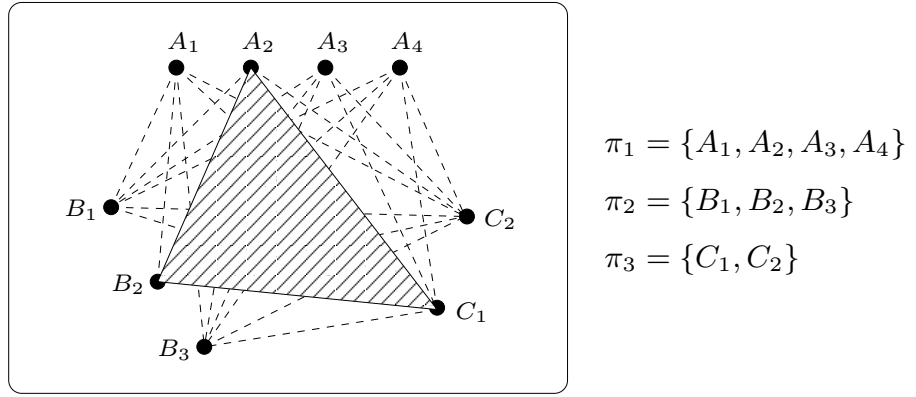
\begin{figure}[ht]
  \begin{center}
    \scalebox{1.3}{\input{rec3_detail_alt.pstex_t}}
  \end{center}
  \caption[\ \ A $3$-rectangle.]{\label{f-rec3_detail_alt}
  A $3$-rectangle $(\pi_1, \pi_2, \pi_3)$,
represented by a $3$-partite graph, and labeled on the triangles.}
\end{figure}

We denote by $\mc{N}$ the species of rectangles, and by
$\mc{N}^{(k)}$ the species of $k$-rectangles.

Let $n=\prod_{i=1}^k {a_i}$, and
let $\Delta$ be the set of bijections of the form
$$\delta:[a_1] \times [a_2]\times\cdots \times [a_k]  \rightarrow [n].$$
Note that the cardinality of the set $\Delta$ is $n!$. The group
$$\prod_{i=1}^k \mathfrak{S}_{a_i}=
\{ \sigma=(\sigma_1, \sigma_2, \dots, \sigma_k):
\sigma_i \in \mathfrak{S}_{a_i}\}$$ acts on the set $\Delta$ by
setting
$$(\sigma\cdot \delta)(i_1, i_2, \dots, i_k)=\delta(\sigma_1(i_1), \sigma_2(i_2),
\dots, \sigma_k(i_k)),$$
for each $(i_1, i_2,\dots, i_k) \in [a_1] \times [a_2]\times\cdots \times [a_k]$.
We observe that this group action result in a set of $\prod_{i=1}^k \mathfrak{S}_{a_i}$-orbits,
and that each orbit consists of exactly $a_1!a_2!\cdots a_k!$ elements of $\Delta$.
Observe further that there is a one-to-one correspondence
between the set of $\prod_{i=1}^k \mathfrak{S}_{a_i}$-orbits on the set $\Delta$
and the set of $k$-rectangles of the form $(\pi_1, \dots, \pi_k)$,
where each $\pi_i$ has $a_i$ blocks.
Therefore, the number of such $k$-rectangles is
$$\biggl\{\aatop{n}{a_1, a_2, \dots, a_k} \biggr\}:=\frac{n!}{a_1!a_2!\cdots a_k!}.$$

\begin{dfn}\label{dfn-arith}
  Let $F_1$ and $F_2$ be species of structures  with $F_1[\emptyset]=F_2[\emptyset]=\emptyset.$
  The \emph{arithmetic
  product} of $F_1$ and $F_2$, denoted $F_1\boxdot F_2$, is
  defined by setting for each finite set $U$,
  $$(F_1\boxdot F_2)[U] = \sum_{(\pi_1\,\pi_2)\in \mc{N}[U]} F_1[\pi_1] \times F_2[\pi_2],$$
  where the sum represents the disjoint union (See Figure~\ref{f-arith_12}).
\begin{figure}[ht]
\begin{center}
\scalebox{.8}{\input{arith_12.pstex_t}}
\end{center}
\caption[\ \ Arithmetic product of species.]{\label{f-arith_12} Arithmetic product $F_1 \boxdot F_2$.}
\end{figure}

  In other words, an $F_1 \boxdot F_2$-structure on a finite set
  $U$ is a tuple of the form $((\pi_1, f_1), (\pi_2, f_2))$, where $(\pi_1, \pi_2)$ is a rectangle on $U$ and $f_i$
  is an $F_i$-structure on the blocks of $\pi_i$ for each $i$.
  A bijection $\sigma: U \rightarrow V$ sends a
  partition $\pi$ of $U$ to a partition $\pi'$ of $V$, namely,
  $\sigma(\pi)=\pi'=\{ \sigma(B) : B \text{ is a block of }
  \pi\}$. Thus $\sigma$ induces a bijection $\sigma_{\pi}:
  \pi
  \rightarrow \pi'$, sending each block of $\pi$ to a block of
  $\pi'$.
  The transport of structures for any bijection $\sigma: U
  \rightarrow V$ is defined by
  \begin{equation}
    (F_1 \boxdot F_2) [\sigma] ((\pi_1, f_1), (\pi_2, f_2)) = ((\pi_1',
    F_1[\sigma_{\pi_1}](f_1)), \,
    (\pi_2', F_2[\sigma_{\pi_2}](f_2))). \notag
  \end{equation}
\end{dfn}

Maia and M\'endez showed that the arithmetic product of species
is commutative, associative, distributive, and with a unit $X$, the species
of singleton sets:
$$F_1 \boxdot X = X \boxdot F_1 = F_1.$$

\begin{dfn}\label{dfn-arithk}
The \emph{arithmetic product} of species  $F_1, F_2, \dots, F_k$
with $F_i(\emptyset)=\emptyset$ for all $i$ is defined by setting
$\mathboxdot_{i=1}^k F_i=F_1 \boxdot F_2 \boxdot \cdots \boxdot
F_k,$ which sends each finite set $U$ to the set
\begin{equation}
  \mathboxdot_{i=1}^k F_i[U]=\sum F_1[\pi_1] \times F_2[\pi_2] \times \cdots \times F_k[\pi_k], \notag
\end{equation}
where the sum is taken over all $k$-rectangles $(\pi_1, \pi_2, \dots, \pi_k)$ of $U$, and
represents the disjoint union. We denote by
$F^{\boxdot k}$ the arithmetic product of $k$ copies of $F$.

For each bijection $\sigma: U \rightarrow V$, the transport of
structures of $\boxdot_{i=1}^k F_i$ along $\sigma$ sends
an $\boxdot_{i=1}^k F_i$-structure on $U$ of the form
$$( (\pi_1, f_1), (\pi_2, f_2), \dots, (\pi_k, f_k) )$$
to an $\boxdot_{i=1}^k F_i$-structure on $V$ of the form
$$( (\pi_1', F_1[\sigma_{\pi_1}]f_1), (\pi_2', F_2[\sigma_{\pi_2}]f_2),
 \dots, (\pi_k', F_k[\sigma_{\pi_k}]f_k) ),$$
where $\sigma_{\pi_i}$ is the bijection induced by $\sigma$
sending blocks of $\pi_i$ to blocks of $\pi_i'$.
\end{dfn}

Maia and M\'endez proved the following proposition which illustrates
that the Dirichlet exponential
generating functions are useful for enumeration involving the arithmetic
product of species.

\begin{prop}
\label{dirichexpsp}
{(Maia and M\'endez)}
Let $F_1$ and $F_2$ be species with
$F_i[\emptyset]=\emptyset$ for $i=1,2$. Then
  \begin{equation}
    \mathfrak{D}({F_1 \boxdot F_2})  = \mathfrak{D}({F_1})\,
    \mathfrak{D}({F_2}).\label{degfarithprod}
  \end{equation}
\end{prop}

\begin{thm}\label{cycleindex-arith}
\emph{(Maia and M\'endez)}
  Let species $F_1$ and $F_2$ satisfy
  $F_1[\emptyset]=F_2[\emptyset]=\emptyset$. Then we have
  \begin{equation}
    Z_{F_1 \boxdot F_2} = Z_{F_1} \boxtimes Z_{F_2}, \label{ciarithprod}
  \end{equation}
  where the operation $\boxtimes$
   on the right-hand side of the
  equation is a bilinear operation on symmetric functions defined by setting
 $$ p_{\nu}:=p_{\lambda} \boxtimes p_{\mu},$$
where
$$c_k(\nu)=\sum_{\lcm(i,j)=k} \gcd(i,j)\, c_i(\lambda) c_j(\mu),$$
in which $\lcm(i,j)$ denotes the \emph{least common multiple} of $i$
and $j$, and $\gcd(i,j)$ denotes the \emph{greatest common
divisor} of $i$ and $j$.
\end{thm}

Furthermore, the arithmetic product of molecular species and the Cartesian product of graphs
are closely related, as shown in the following proposition.

\begin{prop}\label{graphprod}
 Let $G_1$ and $G_2$ be two graphs that are relatively prime to each
other. Then the species associated to the Cartesian product of
$G_1$ and $G_2$ is equivalent to the
arithmetic product of the species associated to $G_1$ and the
species associated to $G_2$. That is,
  \begin{equation}
   \mc{O}_{G_1\odot G_2} =\mc{O}_{G_1} \boxdot \mc{O}_{G_2}
   \label{g1g2prod}
  \end{equation}
\end{prop}

\begin{proof}
  Let $l(G_1)=m$ and $l(G_2)=n$. Then $l(G_1\odot G_2)=mn$.

   Since $G_1$ and $G_2$ are relatively prime, we get
  $$\aut (G_1 \odot G_2)= \aut(G_1) \times \aut (G_2).$$

  Therefore,
  \begin{align*}
    \mc{O}_{G_1\odot G_2}&=\frac{X^{l(G_1 \odot G_2)}}{\aut(G_1\odot G_2)}=\frac{X^{mn}}{\aut (G_1) \times \aut(G_2)}\\
    &=\frac{X^m}{\aut(G_1)} \boxdot \frac{X^n}{\aut(G_2)}=\mc{O}_{G_1} \boxdot \mc{O}_{G_2}.
  \end{align*}
\end{proof}

Note that if $G_1$ and $G_2$ are not relatively prime to each
other, then the species associated to the Cartesian product of
$G_1$ and $G_2$ is generally different from the arithmetic product
of $\mc{O}_{G_1}$ and $\mc{O}_{G_2}$. This is because the
automorphism group of the product of the graphs is no longer the
product of the automorphism groups of the graphs.

\subsection{Exponential Composition of Species}
\label{exp-comp}

%

Let $A$ be a subgroup of $\mathfrak{S}_m$, and let $B$ be a subgroup of
$\mathfrak{S}_n$. The group $B^A$ defined by Definition~\ref{dfn-wreath} acts on the set
of functions from $[m]$ to $[n]$, and hence can be identified with a
subgroup of $\mathfrak{S}_{n^m}$. This gives rise to a molecular species
$X^{n^m}/B^A$, which is defined to be the \emph{exponential composition} of species.
A more general definition is given in the following.

  Let $F$ be a species of structures with $F[\emptyset]=\emptyset$, let $k$ be
   a positive integer, and let $A$ be a subgroup of $\mathfrak{S}_k$.
Recall that an $F^{\boxdot k}$-structure on a finite set $U$ is a
tuple of the form
$$((\pi_1,  f_1), (\pi_2, f_2), \dots, (\pi_k, f_k)),$$
where $(\pi_1, \pi_2, \dots, \pi_k)$ is a $k$-rectangle on $U$, and each $f_i$
is an $F$-structure on the blocks of $\pi_i$. The group $A$ acts on the set
of $F^{\boxdot k}$-structures by permuting the subscripts of $\pi_i$ and
$f_i$, i.e.,
\begin{equation}
  \alpha ((\pi_1,  f_1),  \dots, (\pi_k, f_k))
  =((\pi_{\alpha(1)},  f_{\alpha(1)}),  \dots,
  (\pi_{\alpha(k)}, f_{\alpha(k)})),
  \notag
\end{equation}
where $\alpha$ is an element of $A$,  $(\pi_{\alpha_1}, \pi_{\alpha_2}, \dots,
\pi_{\alpha_k})$ is a $k$-rectangle on $U$, and each $f_{\alpha_i}$
is an $F$-structure on the blocks of $\pi_{\alpha_i}$.
It is easy to check that this action of $A$ on $F^{\boxdot
k}$-structures is natural, that is, it commutes with
any bijection $\sigma: U \rightarrow
V$. Hence we get a quotient species under this group action.

\begin{dfn}\label{dfn-expcompA}
\emph{(Exponential Composition with a Molecular Species)} Let $F$
be a species with $F[\emptyset]=\emptyset.$ We define the
\emph{exponential composition} of $F$ with the molecular species
$X^k/A$ to be the quotient species, denoted $(X^k/A) \langle F
\rangle$, under the group action described in above. That is,
$$(X^k/A)\langle F \rangle := F^{\boxdot k} /A.$$
\end{dfn}

\begin{thm}
\label{molecular-expcomp}
  Let $A$ and $B$ be subgroups of $\mathfrak{S}_m$ and $\mathfrak{S}_n$, respectively,
  and let $B^A$ be the exponentiation group of $A$ with $B$. Then we have
  \begin{equation}
     \frac{X^m}{A} \,
    \biggl\langle \frac{X^n}{B} \biggr\rangle
    = \frac{X^{n^m}}{B^A}. \notag
  \end{equation}
  As a consequence, we have
    \begin{equation}
    Z_{(X^m/A) \langle X^n/B \rangle} = Z(B^A). \notag
  \end{equation}
\end{thm}

\begin{proof}
  Since the arithmetic product is associative, we have
  $$\biggl( \frac{X^n}{B}\biggr)^{\boxdot m}= \frac{X^N}{B^m} ,$$
  where $B^m$ is the the  product of $m$ copies of $B$,
  acting on the set $$\aatop{\underbrace{[n] \times [n] \times \cdots \times [n]}}{m\  \text{copies}} $$
  piecewisely, and hence viewed as a subgroup of $\mathfrak{S}_{n^m}$.
  Therefore, the set of $(X^n/B)^{\boxdot m}$-structures
  on $[N]$ can be identified with the set of $B^m$-orbits of linear orders on $[N]$.

    The group $A$ acts on these $B^m$-orbits of linear orders by permuting the subscripts. This action
    results in the quotient species
  \begin{align*}
    \frac{X^m}{A} \,
    \biggl\langle \frac{X^n}{B} \biggr\rangle
    &= \biggl.\biggl( \frac{X^n}{B}\biggr)^{\boxdot m}  \biggr/A
    = \biggl.\biggl( \frac{X^N}{B^m} \biggr) \biggr/A.
  \end{align*}

  We observe that an $A$-orbit of $B^m$-orbits of linear orders on $[N]$ admits an
  automorphism group isomorphic to the exponentiation group $B^A$, hence the quotient species
  $(X^N/B^m)/A$ is the same as the molecular species $X^N/B^A$.
  Figure~\ref{f-molecular4} illustrates a group action of $A$ on a
  set of $(X^n/B)^{\boxdot m}$-structures.
  \begin{figure}[ht]
    \begin{center}
      \scalebox{.95}{\input{molecular4.pstex_t}}
    \end{center}
    \caption[\ \ The formula $((X^n/B)^{\boxdot m})/A=X^{n^m}/(B^A)$.]{\label{f-molecular4} $((X^n/B)^{\boxdot m})/A=X^{n^m}/(B^A)$.}
  \end{figure}
\end{proof}

\begin{dfn}
\label{dfn-expcompspk}
 Let $k$ be a positive integer, and $F$ a
species with $F[\emptyset]=\emptyset.$ We define the
\emph{exponential composition} of $F$ of \emph{order $k$} to be
the species $$\mc{E}_k \langle F \rangle = F^{\boxdot
k}/\mathfrak{S}_k.$$ We set $\mc{E}_0\langle F \rangle=X$.
\end{dfn}

\begin{dfn}\label{dfn-expcompsp}
\emph{(Exponential Composition of Species)}
Let $F$ be a species with $F[\emptyset]=F[1]=\emptyset$. We define the
\emph{exponential composition} of $F$, denoted $\mc{E}\langle F
\rangle$, to be the sum of $\mc{E}_k\langle F \rangle$ on all
nonnegative integers $k$, i.e.,
\begin{equation}
 \mc{E}\langle F \rangle = \sum_{k\ge 0} \mc{E}_k\langle
F \rangle .
 \notag
\end{equation}
\end{dfn}

The exponential composition of species has properties listed in
the following theorems. Theorem~\ref{thm-dirichexpcomp} gives a
connection between the exponential composition and the Dirichlet
exponential generating function of species.
Theorem~\ref{thm-expcomp-property} lists further properties of the
exponential composition of the sum of two species.

\begin{thm}
  \label{thm-dirichexpcomp}
  Let $F$ be a species with $F[\emptyset]=F[1]=\emptyset$. Then
  \begin{equation}
    \mathfrak{D}(\mc{E}\langle F
    \rangle)=\exp(\mathfrak{D}(F)).\notag
  \end{equation}
\end{thm}

\begin{proof}
Each $F^{\boxdot k}/\mathfrak{S}_k$-structure on a finite set $U$
is an $\mathfrak{S}_k$-orbit of $F^{\boxdot k}$-structures on $U$,
where the action is taken by permuting the subscripts of the
$F^{\boxdot k}$-structures. We observe that there are $k!$
$F^{\boxdot k}$-structures in each of the $\mathfrak{S}_k$-orbits.
Therefore,
$$\left| \frac{F^{\boxdot k}}{\mathfrak{S}_k}\,[n]\right|
=\frac{\left| F^{\boxdot k}[n]  \right|}{k!},$$ and $$\mathfrak{D}
(\mc{E}_k \langle F\rangle)=\mathfrak{D} (F^{\boxdot
k}/\mathfrak{S}_k) = \frac{\mathfrak{D} (F^{\boxdot
k})}{k!}=\frac{\mathfrak{D}(F)^k}{k!}.$$ It follows that
\begin{align}
  \mathfrak{D}(\mc{E} \langle F\rangle)
  =\mathfrak{D}\biggl( \sum_{k\ge 0} \mc{E}_k \langle
  F\rangle)\biggr)
  = \sum_{k \ge 0} \mathfrak{D}(\mc{E}_k \langle F\rangle)
  = \sum_{k \ge 0} \frac{\mathfrak{D}(F)^k}{k!}
  =\exp(\mathfrak{D}(F)).
\end{align}
\end{proof}

\begin{thm}
{(Properties of the Exponential Composition)}
\label{thm-expcomp-property}
  Let $F_1$ and $F_2$ be species with
$F_1[\emptyset]=F_2[\emptyset]=F_1[1]=F_2[1]=\emptyset$, and let
$k$ be any nonnegative integer. Then
\begin{align}
  \mc{E}_k \langle F_1 + F_2 \rangle &= \sum_{i=0}^k \mc{E}_i
  \langle F_1 \rangle \boxdot \mc{E}_{k-i} \langle F_2 \rangle,
  \notag \\
  \mc{E} \langle F_1 + F_2 \rangle &= \mc{E}\langle F_1
  \rangle \boxdot \mc{E} \langle F_2 \rangle. \label{expcomp-distributivity}
\end{align}
\end{thm}

  We observe that an $(F_1+F_2)^{\boxdot k}$-structure on a finite
  set $U$ is a rectangle on $U$ with each partition in the rectangle
  enriched with either an $F_1$ or an $F_2$-structure. Taking the $\mathfrak{S}_k$-orbits
  of these $(F_1+F_2)^{\boxdot k}$-structures on $U$ means basically
  making every partition of the rectangle ``indistinguishable". Hence
  in each $\mathfrak{S}_k$-orbit, all partitions enriched with an $F_1$-structure are grouped
  together to give an $\mathfrak{S}_{k_1}$-orbit of $F_1^{\boxdot k_1}$-structures,
  and the remaining partitions are grouped together to give an $\mathfrak{S}_{k_2}$-orbit
  of $F_2^{\boxdot k_2}$-structures, where $k_1$ and $k_2$ are nonnegative
  integers whose sum is equal to $k$.

\begin{proof}[Proof of Theorem~\ref{thm-expcomp-property}]
First, we prove that for any
nonnegative integer $k$,
\begin{align*}
\mc{E}_k \langle F_1 + F_2 \rangle &= \sum_{i=0}^k \mc{E}_i
  \langle F_1 \rangle \boxdot \mc{E}_{k-i} \langle F_2 \rangle.
\end{align*}

  The case when $k=0$ is trivial. Let us consider $k$ to be a
  positive integer. Let $s$ and $t$ be nonnegative integers whose sum
  equals $k$. Let $U$ be a finite set. To get an $\mc{E}_s \langle F_1 \rangle \boxdot
  \mc{E}_t \langle F_2 \rangle$-structure on $U$, we first
  take a rectangle $(\rho, \tau)$ on $U$, and then take
  an ordered pair $(a, b)$, where $a$ is
  an $\mc{E}_s \langle F_1 \rangle$-structure on the blocks of $\rho$,
  and $b$ is an $\mc{E}_t \langle F_2 \rangle$-structure on the blocks of $\tau$.
  That is,
  \begin{align*}
    a &= \{ (\rho_1, f_1), \dots, (\rho_s, f_s)\}, \qquad
    b = \{ (\tau_1, g_1), \dots, (\tau_t, g_t)\},
  \end{align*}
  where $(\rho_1, \dots, \rho_s)$ is a rectangle on the blocks of $\rho$,
  $(\tau_1, \dots, \tau_t)$ is a rectangle on the blocks of $\tau$, $f_i$
  is an $F_1$-structure on the blocks of $\rho_i$, and $g_j$ is an
  $F_2$-structure on the blocks of $\tau_j$.

  As pointed out by Maia and M\'endez~\cite{arith}, for any nonnegative integers $i,j$,
the species of $(i+j)$-rectangles is isomorphic to the arithmetic product of
the species of $i$-rectangles and the species of $j$-rectangles:
  $$  \mc{N}^{(i+j)} = \mc{N}^{(i)} \boxdot \mc{N}^{(j)}.$$
  It follows that $(\rho_1, \dots, \rho_s, \tau_1, \dots, \tau_t)$
  is a rectangle on $U$.

  On the other hand, let $x$ be an $\mc{E}_k\langle F_1 + F_2\rangle$-structure on
  $U$. We can write $x$ as a set of the form
 $$x=\{(\pi_1,  f_1), \dots, (\pi_r, f_r), (\pi_{r+1}, g_{r+1}), \dots, (\pi_k, g_k)\},$$
 where $(\pi_1, \pi_2, \dots, \pi_k)$ is a $k$-rectangle on $U$,
 $r$ is a nonnegative integer between $0$ and $k$,
 each $f_i$ is an $F_1$-structure on $\pi_i$ for $i=1,\dots, r$, and
 each $g_j$ is an $F_2$-structure on $\pi_j$ for $j=r+1, \dots, k$.

 We then write $x=(x_1, x_2)$, where
 \begin{align*}
   x_1 &= \{(\pi_1,  f_1), \dots, (\pi_r, f_r)\}, \qquad
   x_2 =\{(\pi_{r+1}, g_{r+1}), \dots, (\pi_k, g_k)\}.
 \end{align*}

  Hence running through values of $s$ and $t$, we get that
  the set of  $\mc{E}_s \langle F_1 \rangle \boxdot
  \mc{E}_t \langle F_2 \rangle$-structures on $U$, written in the form of the pairs
  $(a,b)$  whose construction we described
  in above, corresponds naturally to the set of $\mc{E}_k\langle F_1 + F_2\rangle$-structures
  on $U$.

%
%
%
%
%
%
%
The proof of
\begin{equation}
  \mc{E} \langle F_1 + F_2 \rangle = \mc{E}\langle F_1
  \rangle \boxdot \mc{E} \langle F_2 \rangle. \notag
\end{equation}
is straightforward using the properties of the arithmetic product, namely, the
commutativity, associativity and distributivity:
  \begin{align*}
    \mc{E} \langle F_1 + F_2 \rangle
    &= \sum_{k\ge0} \mc{E}_k \langle F_1 + F_2 \rangle
    = \sum_{k\ge0} \sum_{\aatop{i+j=k}{i,j\ge0}} \mc{E}_i \langle F_1 \rangle
       \boxdot \mc{E}_j \langle F_2 \rangle\\
    &= \biggl( \sum_{i\ge0} \mc{E}_i \langle F_1 \rangle \biggr)
       \boxdot \biggl( \sum_{j\ge0} \mc{E}_j \langle F_2 \rangle
       \biggr)= \mc{E}\langle F_1
  \rangle \boxdot \mc{E} \langle F_2 \rangle.
  \end{align*}
\end{proof}

Note that identity~\eqref{expcomp-distributivity} is analogous to the identity
about the composition of a sum of species with the species of sets $\mc{E}$:
$$\mc{E}(F_1+F_2) = \mc{E}(F_1) \, \mc{E}(F_2).$$
What is more, \eqref{expcomp-distributivity} illustrates a kind of
distributivity of the exponential composition. In fact, if a species of
structures $F$ has its molecular decomposition written in the form
\begin{equation}
  F= \sum_{\aatop{M \subseteq F}{M \text{ molecular}}} M, \notag
\end{equation}
then the exponential composition of $F$ can be written as
\begin{equation}
  \expcomp{F}= \mathboxdot_{\aatop{M \subseteq F}{M \text{ molecular}}}\expcomp{M}. \notag
\end{equation}

\subsection{Cycle Index of Exponential Composition}

The cycle index polynomial of the exponentiation group was given by Palmer and
Robinson~\cite{palmer}. They defined the following operators $I_k$ for positive integers $k$.

Let $\mathfrak{R}=\mb{Q}\,[p_1, p_2,\dots]$ be the ring of polynomials
with the operation $\boxtimes$ as defined in
Theorem~\ref{cycleindex-arith}. Palmer and Robinson defined for
positive integers $k$ the $\mb{Q}$-linear operators $I_k$ on
$\mathfrak{R}$ as follows:

Let $\lambda=(\lambda_1, \lambda_2, \dots)$ be a partition of $n$.
The action of $I_k$ on the monomial $p_\lambda$ is given by
\begin{equation}
  I_k(p_\lambda)=p_\gamma, \label{operatorik}
\end{equation}
where $\gamma=(\gamma_1, \gamma_2, \dots)$ is the partition of
$n^k$ with
  \begin{equation}
    c_j(\gamma) =\frac{1}{j} \sum_{l|j} \mu \biggl( \frac{j }{l} \biggr)
    \biggl(\sum_{i\, |\, l/\gcd(k,l)} i c_i(\lambda) \biggr)^{ \gcd(k, l)}.
    \notag
  \end{equation}

Furthermore, $\{ I_k \}$ generates a $\mb{Q}$-algebra $\Omega$ of
$\mb{Q}$-linear operators on $\mathfrak{R}$. For any elements $I, J
\in \Omega$, any $r \in \mathfrak{R}$ and $a \in \mb{Q}$, we set
\begin{align}
(aI)(r)&=a(I (r)), \notag\\
(I+J)(r)&= I(r)+J(r), \notag\\
(IJ)(r)&= I(r) \boxtimes J(r). \label{operatorij}
\end{align}

As discussed in Palmer and Robinson's paper~\cite{palmer}, if $I_m(p_\mu)=p_\nu$, then $\nu$ is the cycle type
of an element $(\alpha, \tau)$ of the exponentiation group $B^A$ acting on
$[n]^m$, where $\alpha$ is a permutation in $A$ with
a single $m$-cycle, and $\tau \in B^m$ is such that
$\mu$ is the cycle type of the permutation $\tau(m)\tau(m-1)\cdots \tau(2)\tau(1)$.

\begin{dfn}\label{dfn-expcompcycle}
Let $f_1$ and $f_2$ be elements of the ring $\mathfrak{R}=\mb{Q}\,[p_1, p_2,\dots]$.
We define the \emph{exponential composition} of
$f_1$ and $f_2$, denoted $f_1 \ast f_2$, to be the
  image of $f_2$ under the operator obtained by substituting the
  operator $I_r$ for the variables $p_r$ in $f_1$.
\end{dfn}

Note that the operation $\ast$ is linear in the left parameters, but not on the right parameters.
We call this the \emph{partial linearity} of the operation $\ast$.

Let $A$ be a subgroup of $\mathfrak{S}_m$, and let $B$ be a subgroup
of $\mathfrak{S}_n$. Palmer and Robinson~\cite[pp.~128--131]{palmer}
proved that the cycle index polynomial of $B^A$ is
  the exponential composition of $Z(A)$ with $Z(B)$. That is,
  \begin{equation}
    Z(B^A)=Z(A) \ast Z(B).
    \notag
  \end{equation}

As a consequence of Theorem~\ref{molecular-expcomp},
we get the cycle index of the species $(X^m/A) \langle X^n/B \rangle$:
\begin{equation}
 Z_{(X^m/A) \langle X^n/B \rangle} = Z(A) \ast Z(B). \notag
\end{equation}

Next we generalize Palmer and Robinson's result to get the formula for the cycle index
of the exponential composition of an arbitrary species. First, we introduce a lemma
that is a generalization of the Cauchy-Frobenius Theorem, alias
Burnside's Lemma. For the proof of a more general result, with
applications and further references, see Robinson~\cite{robinson}.
Another application is given in~\cite{functorial}.

\begin{lem}\label{burnsides}
\emph{(Cauchy-Frobenius)}
  Suppose that a finite group $M \times N$ acts on a set $S$. The groups $M$ and
  $N$, considered as subgroups of $M \times N$, also act on $S$.
   The group $N$ acts on the set of $M$-orbits.
  Then for any $g \in N$,
  the number of $M$-orbits fixed by $g$ is given by
  \begin{equation}
   \frac {1}{|M|} \, \sum_{f \in M} \fix(f,g), \notag
  \end{equation}
  where $\fix(f,g)$ denotes the number of elements in $S$ that are
  fixed by $(f,g)\in M \times N$.
\end{lem}

\begin{thm}\label{expAF}
\emph{(Cycle Index of the Exponential Composition)}
  Let $A$ be a subgroup of $\mathfrak{S}_k$, and let $F$ be a species of
  structures concentrated on the cardinality $n$. Then the cycle
  index of the species $(X^k/A) \langle F \rangle$ is given by
  \begin{equation}
    Z_{(X^k/A) \langle F \rangle}=Z(A) \ast Z_F,
    \label{expecompspeciesformula0}
  \end{equation}
   where the expression $Z(A) \ast Z_F$ denotes
  the  image of $Z_F$ under
  the operator obtained by substituting the
  operator $I_r$ for the variables $p_r$ in $Z(A)$.
\end{thm}

\begin{rem}
[\emph{Notation and Set-up}]
\label{rem-thmAF}
We denote by $\Par_n$ the set of partitions of $n$, and by $\Par_n^k$
the set of $k$-sequences of partitions of $n$.

For  fixed integers $n$, $k$, and $N=n^k$,
we  denote by $\mc{N}_N$
the species of \emph{$k$-dimensional cubes, or $k$-cubes, on $[N]$}, defined by
$$\mc{N}_N = \mc{E}_n^{\boxdot k}[N].$$
We also call the elements of the set $(X^n)^{\boxdot k}[N]$ \emph{$k$-dimensional ordered cubes on $[N]$}.

Let $\sigma$ be a permutation on $[k]$ with cycle type
$$\c.t.(\sigma)=(r_1,  r_2 \dots, r_d).$$ Then $\sigma$ acts on
the $F^{\boxdot k}$-structures by permuting
  the subscripts. Let $\nu$ be a partition of $N$. Let $\delta$ be a permutation
of $[N]$ with cycle type $\nu$. Then $\delta$ acts on the
$F^{\boxdot k}$-structures by transport of structures.
We also introduce the notation
$$I(\c.t.(\sigma); \lambda^{(1)}, \lambda^{(2)}, \dots, \lambda^{(d)})=
I_{r_1}(p_{\lambda^{(1)}})\boxtimes I_{r_2}(p_{\lambda^{(2)}}) \boxtimes
\cdots \boxtimes I_{r_d}(p_{\lambda^{(d)}}).$$
We denote by $\Rec_F(\sigma,\nu)$
a function on the pair $(\sigma, \nu)$ defined by
  \begin{equation}
  \Rec_F(\sigma, \nu) := \sum\,\frac{\prod_{i=1}^d \fix\, F[\lambda^{(i)}]}{z_{ \lambda^{(1)}}\cdots z_{ \lambda^{(d)}}}, \label{dfn-recf}
    \end{equation}
where the summation is over all sequences $(\lambda^{(1)}, \lambda^{(2)}, \dots, \lambda^{(d)})$ in $ \Par_n^d$
with $$I(\c.t.(\sigma); \lambda^{(1)}, \lambda^{(2)}, \dots, \lambda^{(d)}))= p_\nu.$$
We denote by $\fix_F(\sigma, \delta)$  the number of $F^{\boxdot k}$-structures on the set $[N]$ fixed by the
  joint action of the pair $(\sigma, \delta)$.

\end{rem}

\begin{proof}[Proof of Theorem~\ref{expAF}]

Let $\nu$ be a partition of $N$. It suffices to prove that the coefficients
of $p_\nu$ on both sides of Equation~\eqref{expecompspeciesformula0}
are equal.

The right-hand side of Equation~\eqref{expecompspeciesformula0} is
    \begin{align*}
    \ Z(A) \ast Z_F
    &= \biggl( \frac{1}{|A|}\, \sum_{\sigma \in A} p_{\c.t.(\sigma)} \biggr)
    \ast \biggl( \sum_{\lambda \vdash n}\, \fix\, F[\lambda]\,\frac{p_\lambda}{z_\lambda} \biggr)\\
    &= \frac{1}{|A|}\, \sum_{\sigma \in A}
    I_{\c.t.(\sigma)} \biggl( \sum_{\lambda \vdash n}\, \fix\, F[\lambda]\,\frac{p_\lambda}{z_\lambda} \biggr).
    \end{align*}
    For $\sigma \in A$ with $\c.t.(\sigma)=(r_1, r_2, \dots, r_d)$, we have
    $$I_{\c.t.(\sigma)}=I_{r_1}\cdots I_{r_d},$$ and
    \begin{multline*}
     I_{\c.t.(\sigma)} \biggl( \sum_{\lambda \vdash n}\,
    \fix\, F[\lambda]\,\frac{p_\lambda}{z_\lambda} \biggr)\\
    = I_{r_1}\biggl( \sum_{\lambda \vdash n}\, \fix\, F[\lambda]\,\frac{p_\lambda}{z_\lambda} \biggr)
     \boxtimes
     I_{r_2}\biggl( \sum_{\lambda \vdash n}\, \fix\, F[\lambda]\,\frac{p_\lambda}{z_\lambda} \biggr)
     \boxtimes
    \cdots \boxtimes I_{r_d}\biggl( \sum_{\lambda \vdash n}\, \fix\, F[\lambda]\,\frac{p_\lambda}{z_\lambda} \biggr).
    \end{multline*}

        Therefore, the coefficient of $p_\nu$ in the expression
    $Z(A) \ast Z_F$ is
    \begin{align}
    & \frac{1}{|A|}\, \sum \,
    \frac{\prod_{i=1}^d \fix F[\lambda^{(i)}]}{z_{\lambda^{(1)}}\cdots z_{\lambda^{(d)}}}
    =\frac{1}{|A|}\, \sum_{\sigma \in A} \Rec_F(\sigma, \nu), \label{thmAF-RHS}
  \end{align}
  where the summation on the left-hand side is taken over
 all sequences $(\lambda^{(1)}, \lambda^{(2)}, \dots, \lambda^{(d)})$ in $ \Par_n^d$ for some $d \ge 1$
 and all $\sigma \in A$ with $\c.t.(\sigma)=(r_1, r_2, \dots, r_d)$ such that
    $$I(\c.t.(\sigma); \lambda^{(1)}, \lambda^{(2)}, \dots, \lambda^{(d)})= p_\nu,$$
  and $\Rec_F(\sigma, \nu)$ on the right-hand side is as defined by~\eqref{dfn-recf} in Remark~\ref{rem-thmAF}.

The left-hand side of Equation~\eqref{expecompspeciesformula0} is
  \begin{align*}
    Z_{F^{\boxdot k}/A}& =\sum_{\nu \vdash N} \fix\,\frac{F^{\boxdot k}}{A}\,[\nu]\, \frac{p_\nu}{z_\nu}.
  \end{align*}
  Therefore, the coefficient of $p_\nu$ in the expression
  $Z_{F^{\boxdot k}/A}$ is
  \begin{align*}
   & \frac{1}{z_\nu} \, \fix\,\frac{F^{\boxdot k}}{A}\,[\nu].
  \end{align*}
  We then apply Theorem~\ref{burnsides} to get that the number of $A$-orbits of
  $F^{\boxdot k}$-structures on $[N]$ fixed by a permutation $\delta \in
  \mathfrak{S}_N$ of cycle type $\nu$ is
  \begin{align}
    \fix\,\frac{F^{\boxdot k}}{A}\,[\nu] &= \fix\,\frac{F^{\boxdot k}}{A}\,[\delta]=
    \frac{1}{|A|} \, \sum_{\sigma \in A} \fix_F(\sigma, \delta), \label{thmAF-LHS}
  \end{align}
  where $\fix_F(\sigma,\delta)$ is as defined in Remark~\ref{rem-thmAF}.

  Therefore, combining~\eqref{thmAF-LHS} and~\eqref{thmAF-RHS},
  the proof of Equation~\eqref{expecompspeciesformula0} is reduced
  to showing that
  \begin{equation}
    \fix_F(\sigma,\delta)= z_\nu \Rec_F(\sigma, \nu), \label{proof-recfix}
  \end{equation}
  for any $\delta, \nu$ and $\sigma$.

To prove~\eqref{proof-recfix},
we start with observing that in order for an $F^{\boxdot k}$-structure on $[N]$
of the form $$((\pi_1, f_1), (\pi_2, f_2), \dots, (\pi_k, f_k) )$$ to be
fixed by the pair $(\sigma, \delta)$,
it is necessary that $(\sigma, \delta)$ fixes the $k$-cube
of the form $(\pi_1, \pi_2, \dots, \pi_k) \in \mc{N}_N$.
This is equivalent to saying that
  \begin{equation}
  \mc{N}_N[\delta] (\pi_1, \pi_2, \dots, \pi_k)=
  (\pi_{\sigma (1)}, \pi_{\sigma (2)}, \dots, \pi_{\sigma (k)}). \label{proof-sigmadelta}
  \end{equation}

Suppose~\eqref{proof-sigmadelta} holds for some $k$-cube $(\pi_1, \pi_2, \dots, \pi_k) \in \mc{N}_N$.
We let $\beta_i \in \mathfrak{S}_n$ be the induced action of $\delta$ on the blocks of $\pi_i$, for $i=1, 2, \dots, k$.
That is, $$\mc{N}_N[\delta](\pi_i)= \beta_i (\pi_{\sigma(i)})$$ for all $i \in [k]$.
\begin{figure}[ht]
  \begin{center}
    \scalebox{.8}{\input{pi_sigma_beta.pstex_t}}
  \end{center}
\end{figure}

Now we consider the simpler case when  $\sigma$ is a $k$-cycle, say, $\sigma=(1,2,\dots, k)$. Then
 the action
of $\delta$ sends $(\pi_1, \pi_2, \dots, \pi_k)$ to $(\pi_2, \pi_3, \dots, \pi_1)$.
Let $\beta=\beta_1\beta_2 \cdots \beta_k$. The above discussion is saying that
$$I_k (p_{\c.t.(\beta)})=p_\nu.$$

On the other hand, given a partition $\lambda$ of $n$ satisfying $I_k (p_\lambda)=p_\nu$,
there are $n!/z_\lambda$ permutations in $\mathfrak{S}_n$
with cycle type $\lambda$. Let $\beta$ be one of such. Then
the number of sequences $(\beta_1, \beta_2, \dots, \beta_k)$ whose product equals
$\beta$ is $(n!)^{k-1}$, since we can choose $\beta_1$ up to $\beta_{k-1}$ freely,
and $\beta_k$ is thereforee determined.
All such sequences $(\beta_1, \beta_2, \dots,\beta_k)$
will satisfy $I_k (p_{\c.t.(\beta_1\cdots\beta_k)})=p_\nu$,
thus their action on an arbitrary $k$-dimensional ordered cube, combined with the action
of $\sigma$ on the subscripts, would result in a permutation on $[N]$ with cycle type $\nu$.
But there are $N!/z_\nu$ permutations with cycle type $\nu$, and only one of them is
the $\delta$ that we started with. Considering that
the $k$-cubes are just $\mathfrak{S}_n^k$-orbits of the $k$-dimensional ordered
cubes, we count the number of $k$-cubes that are fixed by the pair $(\sigma, \delta)$
with the further condition that the product of the induced permutations on the $\pi_i$ by $\delta$
has cycle type $\lambda$:
\begin{multline*}
\qquad \quad
  \frac{\# \big\{ \aatop{ (\beta_1, \beta_2, \dots, \beta_k) \in \mathfrak{S}_n^k }{\text{ with } \c.t.(\beta_1 \cdots \beta_k)
  =\lambda  }\big\} \cdot \# \big\{ k\text{-dimensional ordered cubes} \big\}}
  {\# \big\{ \aatop{\text{ permutations on }[N] }{\text{ with cycle type }\nu}  \big\}
  \cdot\# \big\{ \aatop{k\text{-dimensional ordered cubes }}{\text{in each equivalence class}}  \big\}}\\
 =  \frac{[(n!)^{k-1}\cdot n!/z_\lambda] \cdot N!}{N!/z_\nu \cdot (n!)^k}
 =\frac{z_\nu}{z_\lambda}.
 \qquad \quad
\end{multline*}

Now we try to compute how many $F^{\boxdot k}$-structures of the form
$$((\pi_1, f_1), (\pi_2, f_2),  \dots, (\pi_k, f_k)),$$
based on a given rectangle $(\pi_1, \pi_2, \dots, \pi_k)$ that is fixed by the $\beta_i$
with $$\c.t.\biggl(\prod_i \beta_i \biggr)=\lambda,$$
are fixed by the pair $(\sigma, \delta)$.
We observe that the action of $(\sigma,\delta)$ determines
that $f_k=F[\beta_1]f_1$ and $f_i=F[\beta_{i-1}]f_{i-1}$ for $i=2,3,\dots,k$, and hence
$$f_k=F[\beta_1]F[\beta_2]\cdots F[\beta_k] f_k=F[\beta]f_k=F[\lambda]f_k.$$
In other words, $$f_k \in \Fix F[\lambda].$$
Hence as long as we choose an $f_k$ from $\Fix F[\lambda]$, then all the other
$f_i$ for $i<k$ are determined by our choice of $f_k$. There are $\fix F[\lambda]$
such choices for $f_k$.

Therefore, in the case when $\sigma$ is a $k$-cycle,
 we get that the number of  $F^{\boxdot k}$-structures on the set $[N]$
fixed by the pair $(\sigma, \delta)$ is
$$\fix_F(\sigma, \delta)=\sum_{\aatop{\lambda \vdash n}
{I_k(p_\lambda)=p_\nu}} \fix\, F[\lambda] \frac{z_\nu}{z_\lambda}
=z_\nu\Rec_F(\sigma, \nu).$$

Now let us consider the general case when $\sigma$ contains $d$ cycles of
lengths $r_1, r_2, \dots, r_d$.
Let $(\pi_1, \pi_2, \dots, \pi_k)$ be a $k$-cube fixed by the pair $(\sigma, \delta)$.
Again we have~\eqref{proof-sigmadelta}, and we get
an induced $\beta_i$ on the blocks of $\pi_{\sigma^{-1}i}$ for each $i$.

We observe that the action of $\sigma$ on the subscripts of the $k$-cube partitions the list
$\pi_1, \pi_2, \dots, \pi_k$ into $d$ parts, of lengths $r_1, r_2, \dots, r_d$, within each
of which we get a $r_i$-cycle.  We group the
$\beta_i$ on each of the $d$ parts and get $d$ permutations in the group $\mathfrak{S}_n$,
whose cycle types are denoted by $\lambda^{(1)}, \lambda^{(2)}, \dots, \lambda^{(d)}$. This construction
gives that such a sequence of partitions $(\lambda^{(1)}, \lambda^{(2)}, \dots, \lambda^{(d)})$ will
be those that satisfy
$$I(\c.t.(\sigma); \lambda^{(1)}, \lambda^{(2)},  \dots,\lambda^{(d)})=p_\nu.$$

Therefore, the number of $k$-cubes fixed by $(\sigma,\delta)$ corresponding to
such a sequence of partitions $(\lambda^{(1)}, \lambda^{(2)},  \dots, \lambda^{(d)})$ is
\begin{multline*}
  \qquad
  \frac{(n!)^{r_1-1}\cdot n!/z_{\lambda^{(1)}} \cdot \cdots \cdot
  (n!)^{r_d-1}\cdot n!/z_{\lambda^{(d)}} }{N!/z_\nu}\cdot \frac{N!}{(n!)^k}\\
  =  \frac{(n!)^{r_1+ \cdots + r_d}}{N!} \frac{z_\nu}{ z_{\lambda^{(1)}}\cdots z_{\lambda^{(d)}} }=\frac{z_\nu}{ z_{\lambda^{(1)}}\cdots z_{\lambda^{(d)}} }.
  \qquad
\end{multline*}

 The number of $F$-structures that are assigned to this $k$-cube $(\pi_1, \pi_2, \dots, \pi_k)$
that will be fixed under the action of the pair $(\sigma, \delta)$ corresponding to
the sequence of partitions $(\lambda^{(1)}, \lambda^{(2)},  \dots, \lambda^{(d)})$ is hence
$$\fix F[\lambda^{(1)}]\cdots \fix F[\lambda^{(d)}],$$ since, similarly to our previous
discussion, within each of the
$d$ parts, we only need to pick an $F$-structure that is fixed by a permutation of cycle type
$\lambda^{(i)}$, and all other $F$-structures are left determined.

Therefore, we get that for any pair $(\sigma, \delta)$,
  \begin{equation}
    \fix_F(\sigma,\delta)
    =z_\nu \Rec_F(\sigma, \nu), \notag
  \end{equation}
which concludes our proof.
\end{proof}

\begin{rem}
  We can use the molecular decomposition to define \emph{the exponential composition
  of a species $F$ with a species $H$}.
  That is, if the molecular decomposition of $H$ is given by
  $$H=\sum_{\aatop{M \subseteq H}{M \text{ molecular}}} M,$$
  then we define $H \langle F \rangle$ by
  $$H \langle F \rangle =\sum_{\aatop{M \subseteq H}{M \text{ molecular}}} M \langle F \rangle.$$
  The left-linearity of the operation $\ast$ gives that the cycle index of $H\langle F \rangle$ is
  \begin{align*}
  Z_{H\langle F \rangle}&= Z_H \ast Z_F
  = \biggl(\sum_{\aatop{M \subseteq H}{M \text{ molecular}}} Z_M \biggr)  \ast Z_F
  = \sum_{\aatop{M \subseteq H}{M \text{ molecular}}} Z_M \ast Z_F.
  \end{align*}
\end{rem}

\subsection{Cycle Index of the Species of Prime Graphs}

Now we are ready to come back to the species of prime graphs.

\begin{lem}\label{primepower}
Let $P$ be any prime graph, and let $k$ be any nonnegative integer. Then
the species associated to the $k$-th power of $P$
is the exponential composition of $\mc{O}_P$ of
order $k$. That is,
\begin{equation}
  \mc{O}_{P^k}=\mc{E}_k \langle \mc{O}_P \rangle. \notag
\end{equation}
\end{lem}

\begin{proof}
We apply Theorem~\ref{molecular-expcomp} and get
\begin{align*}
\mc{E}_k \langle \mc{O}_P \rangle &= \mc{O}_{P}^{\boxdot k} /
\mathfrak{S}_k = \biggl.\biggl( \frac{X^n}{\aut(P)} \biggr)^{\boxdot k} \biggr/\mathfrak{S}_k= \frac{X^{n^k}}{\aut(P)^{\mathfrak{S}_k}}.
\end{align*}
It follows from Proposition~\ref{primepower_aut} that
$$\mc{E}_k \langle \mc{O}_P \rangle = \frac{X^{n^k}}{\aut(P^k)}=\mc{O}_{P^k}.$$
\end{proof}

We can verify Lemma~\ref{primepower} in an intuitive way. Note that
 the set of $\mc{E}_k \langle \mc{O}_P \rangle$-structures on a finite set
$U$ is the set of $\mathfrak{S}_k$-orbits of  $\mc{O}_P^{\boxdot k}$-structures
on $U$, and an element of  $\mc{E}_k \langle \mc{O}_P \rangle[U]$ of the form
$\{(\pi_1, f_1), \dots, (\pi_k, f_k)\}$
is such that $(\pi_1, \pi_2, \dots, \pi_k)$ is a $k$-rectangle on $U$, and each $f_i$ is
a graph isomorphic to $P$ whose vertex set equal to the blocks of
$\pi_i$. Such a set  $\{(\pi_1, f_1), \dots, (\pi_k, f_k)\}$
corresponds to a graph $G$ isomorphic  to
$P^k$ with vertex set $U$. More precisely, $G$ is  the Cartesian product of the $f_i$
in which each vertex $u \in U$ is of the form $u=B_1 \cap B_2 \cap \cdots \cap B_k$,
where each $B_i$ is one of the
blocks of $\pi_i$.  In this way, we get a one-to-one correspondence between
the $\mc{E}_k \langle \mc{O}_P \rangle$-structures on $U$ and the set of graphs isomorphic to $P^k$
with vertex set $U$.

\begin{thm}
\label{thm-cycleprime}
The species ${\mc{G}^c}$ of connected graphs and ${\mc{P}}$ of prime graphs
satisfy
\begin{equation}
{\mc{G}^c} = \mc{E} \langle {\mc{P}} \rangle. \notag
\end{equation}
\end{thm}

\begin{proof}
In this proof, all graphs considered are unlabeled.

The molecular decomposition
of the species of prime graphs
is $${\mc{P}}=\sum_{P  \text{ prime}} \mc{O}_P,$$
  where each $\mc{O}_P$ is a molecular species which is isomorphic to $X^{l(P)}/\aut(P).$

  Let
  $\{ P_1, P_2, \dots \}$ be the set of unlabeled prime graphs. We have
  \begin{align*}
    \expcomp{{\mc{P}}}
    &= \mc{E} \langle \mc{O}_{P_1} + \mc{O}_{P_2} + \cdots\rangle\\
    &= \expcomp{\mc{O}_{P_1}} \boxdot \expcomp{\mc{O}_{P_2}} \boxdot \cdots \\
    &= (X + \mc{O}_{P_1} + \mc{O}_{P_1^2} + \cdots )\boxdot (X + \mc{O}_{P_2} + \mc{O}_{P_2^2} + \cdots )
    \boxdot \cdots\\
    &= \sum_{i_1, i_2, \dots \ge 0} \mc{O}_{P_1^{i_1}} \boxdot \mc{O}_{P_2^{i_2}} \boxdot \cdots\\
    &= \sum_{i_1, i_2, \dots \ge 0} \mc{O}_{P_1^{i_1} \boxdot P_2^{i_2} \boxdot \cdots } \\
    &= \sum_{C {\text{ connected}}} \mc{O}_C \\
    &= {\mc{G}^c}.
  \end{align*}
\end{proof}

  Note that Theorem~\ref{thm-degfconnprime} follows as a corollary
  of Theorems~\ref{thm-dirichexpcomp} and~\ref{thm-cycleprime}.

\begin{rem}
Recall that the exponential composition of a species $F$ is the
sum of $\mc{E}_k \langle F \rangle$ on all nonnegative integers
$k$:
$$\mc{E}\langle F\rangle=\mc{E}_0
\langle F \rangle + \mc{E}_1 \langle F \rangle +\mc{E}_2\langle
F\rangle+\cdots =X+F+\mc{E}_2\langle F\rangle+\cdots.$$
Theorem~\ref{thm-cycleprime} gives that
\begin{align*}
\mc{G}^c &= X + \mc{P}+\text{higher terms},\\
\mc{P}&=\mc{G}^c -X - \text{higher terms},\\
Z_{\mc{P}}&=Z_{\mc{G}^c} - p_1 - \text{higher terms}.
\end{align*}
Now we can compute the cycle index of the species of prime graphs
$Z_{\mc{P}}$ from the cycle index of the species of connected
graphs $Z_{\mc{G}^c}$, given by formula~\eqref{eq-cycle-conn},
recursively using maple:
\begin{align*}
    Z_{{\mc{P}}} &= \biggl( \frac{1}{2}\,p_1^2 + \,\frac{1}{2}\,p_2  \biggr) +
     \biggl( \frac{2}{3}\,p_1^3 + p_1p_2+\,\frac{1}{3}\,p_3 \biggr) \\
    & \phantom{=}\ \,  + \biggl( \frac{35}{24}p_1^4 + \,\frac{7}{4}\,p_1^2p_2
    +\,\frac{2}{3}\,p_1p_3 +\,\frac{7}{8}\,p_2^2+\,\frac{1}{4}\,p_4  \biggr) \\
    & \phantom{=}\ \,  + \biggl( \frac{91}{15}\,p_1^5 +\,\frac{19}{3}\,p_1^3p_2
    +\,\frac{4}{3}\,p_1^3p_3+5p_1p_2^2+p_1p_4+\,\frac{2}{3}\,p_2p_3+\,\frac{3}{5}\,p_5 \biggr) \\
    & \phantom{=}\ \,  + \biggl( \frac{1654}{45}\,p_1^6 + \,\frac{91}{3}\,p_1^4p_2
    +\,\frac{38}{9}\,p_1^3p_3 +21p_1^2p_2^2 + 2p_1^2p_4 +\,\frac{8}{3}\,p_1p_2p_3
      \biggr.\\
    &  \phantom{=}\ \,  + \biggl.  \frac{4}{5}\,p_1p_5 +\,\frac{47}{6}\,p_2^3
    +\,\frac{5}{2}\,p_2p_4+\,\frac{11}{9}\,p_3^2+\,\frac{2}{3}\,p_6   \biggr) + \cdots.
\end{align*}
\end{rem}

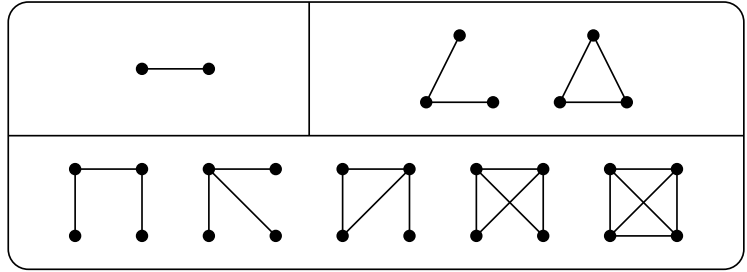
\begin{figure}[ht]
  \begin{center}
    \scalebox{1.4}{\input{molecular_prime.pstex_t}}
  \end{center}
  \caption[\ \ Unlabeled prime graphs.]{\label{f-molecular_prime}
  Unlabeled prime graphs with $n$ vertices, $n \le 4$.}
\end{figure}

Figure~\ref{f-molecular_prime} shows the unlabeled prime graphs on no more than $4$ vertices.

Hence we write down the beginning terms of the molecular decomposition of the species
${\mc{P}}$:
\begin{align*}
  {\mc{P}} & = \mc{E}_2 + (X\cdot \mc{E}_2 + \mc{E}_3)
  + (\mc{E}_2 \circ X^2 + X\cdot \mc{E}_3 + X^2 \cdot \mc{E}_2 + \mc{E}_2 \cdot \mc{E}_2
  + \mc{E}_4) +\cdots.
\end{align*}

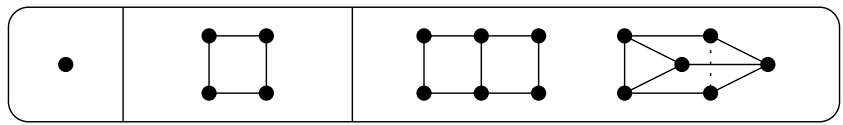
\begin{figure}[ht]
  \begin{center}
    \scalebox{1.2}{\input{non-prime.pstex_t}}
  \end{center}
  \caption[\ \ Unlabeled non-prime graphs.]{\label{f-non-prime} Unlabeled non-prime graphs on $n$ vertices, $n\le 6$.}
\end{figure}

Comparing Figure~\ref{f-molecular_prime} with unlabeled connected graphs with no
more than $4$ vertices,
we see that there is only one unlabeled connected
graph with $4$ vertices that is not prime.
In fact, if we compare the first several terms
of $Z_{\mc{G}^c}$, given by \eqref{eq-cycle-conn}, and $Z_{{\mc{P}}}$ of order no more than $6$, we get that
\begin{align*}
  Z_{\mc{G}^c} - Z_{{\mc{P}}} &=
  p_1 +  \,\frac{1}{8}\,\biggl( p_1^4 + 2p_1^2p_2 + 3\,p_2^2+
  2p_4 \biggr) \\
  & \phantom{=} \ \, + \,\frac{1}{4}\,\biggl( p_1^6 + p_1^2p_2^2+2p_2^3\biggr) +
  \,\frac{1}{12}\,\biggl( p_1^6 + 3p_1^2p_2^2 + 4p_2^3 +2p_3^2+2p_6  \biggr)\cdots,
\end{align*}
which is the cycle index of connected non-prime graphs on no more
than $6$ vertices, as shown in Figure~\ref{f-non-prime}, which
consist of a single vertex, a graph with $4$ vertices, and two
graphs with $6$ vertices.

\appendix
\begin{table}[ht]
\centering
\caption{\label{t-lpg}
Numbers of labeled and unlabeled prime graphs on $n$ vertices,
denoted $p^l_n$ and $p^u_n$, respectively, for $n \le 16$.}
\begin{tabular}{rrr}
$n$   & $p^l_n$                              &     $p^u_n$                  \\ \hline
1     &       0                              &      0                       \\
2     &       1                              &      1                       \\
3     &       4                              &      2                       \\
4     &       35                             &      5                       \\
5     &       728                            &      21                      \\
6     &      26464                           &     110                      \\
7     &     1866256                          &    853                       \\
8     &    251518352                         &   11111                      \\
9     &  66296210432                         &  261077                      \\
10    & 34496477587456                       & 11716550                     \\
11    & 35641657548953344                    & 1006700565                   \\
12    & 73354596197458024448                 & 164059830354                 \\
13    & 301272202649664088951808             & 50335907869219               \\
14    & 2471648811030427594714599424         & 29003487462847208            \\
15    & 40527680937730480229320939012096     & 31397381142761241918         \\
16    & 1328578958335783200943054119287117312& 6396956011322517616514       \\
\hline
\end{tabular}
\end{table}

\bibliographystyle{Amsplain}
\bibliography{prime}

\end{document}

%% file: cart_prod_1.pstex_t
\begin{picture}(0,0)%
\includegraphics{cart_prod_1.pstex}%
\end{picture}%
\setlength{\unitlength}{2565sp}%
\begingroup\makeatletter\ifx\SetFigFont\undefined%
\gdef\SetFigFont#1#2#3#4#5{%
  \reset@font\fontsize{#1}{#2pt}%
  \fontfamily{#3}\fontseries{#4}\fontshape{#5}%
  \selectfont}%
\fi\endgroup%
\begin{picture}(7674,3849)(514,-3298)
\put(1153,-2734){\makebox(0,0)[lb]{\smash{{\SetFigFont{7}{8.4}{\rmdefault}{\mddefault}{\updefault}{\color[rgb]{0,0,0}$b$}%
}}}}
\put(2053,-2734){\makebox(0,0)[lb]{\smash{{\SetFigFont{7}{8.4}{\rmdefault}{\mddefault}{\updefault}{\color[rgb]{0,0,0}$c$}%
}}}}
\put(1003,191){\makebox(0,0)[lb]{\smash{{\SetFigFont{7}{8.4}{\rmdefault}{\mddefault}{\updefault}{\color[rgb]{0,0,0}$1$}%
}}}}
\put(853,-859){\makebox(0,0)[lb]{\smash{{\SetFigFont{7}{8.4}{\rmdefault}{\mddefault}{\updefault}{\color[rgb]{0,0,0}$3$}%
}}}}
\put(1303,-259){\makebox(0,0)[lb]{\smash{{\SetFigFont{7}{8.4}{\rmdefault}{\mddefault}{\updefault}{\color[rgb]{0,0,0}$2$}%
}}}}
\put(2278,-259){\makebox(0,0)[lb]{\smash{{\SetFigFont{7}{8.4}{\rmdefault}{\mddefault}{\updefault}{\color[rgb]{0,0,0}$4$}%
}}}}
\put(1603,-2134){\makebox(0,0)[lb]{\smash{{\SetFigFont{7}{8.4}{\rmdefault}{\mddefault}{\updefault}{\color[rgb]{0,0,0}$a$}%
}}}}
\put(4501,-2221){\makebox(0,0)[lb]{\smash{{\SetFigFont{8}{9.6}{\rmdefault}{\mddefault}{\updefault}{\color[rgb]{0,0,0}$2,c$}%
}}}}
\put(5026,-571){\makebox(0,0)[lb]{\smash{{\SetFigFont{8}{9.6}{\rmdefault}{\mddefault}{\updefault}{\color[rgb]{0,0,0}$2,a$}%
}}}}
\put(5926,-2221){\makebox(0,0)[lb]{\smash{{\SetFigFont{8}{9.6}{\rmdefault}{\mddefault}{\updefault}{\color[rgb]{0,0,0}$3,b$}%
}}}}
\put(7426,-1621){\makebox(0,0)[lb]{\smash{{\SetFigFont{8}{9.6}{\rmdefault}{\mddefault}{\updefault}{\color[rgb]{0,0,0}$4,b$}%
}}}}
\put(6301,-1621){\makebox(0,0)[lb]{\smash{{\SetFigFont{8}{9.6}{\rmdefault}{\mddefault}{\updefault}{\color[rgb]{0,0,0}$2,b$}%
}}}}
\put(6076,-571){\makebox(0,0)[lb]{\smash{{\SetFigFont{8}{9.6}{\rmdefault}{\mddefault}{\updefault}{\color[rgb]{0,0,0}$4,a$}%
}}}}
\put(6226,-1096){\makebox(0,0)[lb]{\smash{{\SetFigFont{8}{9.6}{\rmdefault}{\mddefault}{\updefault}{\color[rgb]{0,0,0}$1,b$}%
}}}}
\put(3901,-2671){\makebox(0,0)[lb]{\smash{{\SetFigFont{8}{9.6}{\rmdefault}{\mddefault}{\updefault}{\color[rgb]{0,0,0}$3,c$}%
}}}}
\put(4501,-1096){\makebox(0,0)[lb]{\smash{{\SetFigFont{8}{9.6}{\rmdefault}{\mddefault}{\updefault}{\color[rgb]{0,0,0}$3,a$}%
}}}}
\put(4726,-121){\makebox(0,0)[lb]{\smash{{\SetFigFont{8}{9.6}{\rmdefault}{\mddefault}{\updefault}{\color[rgb]{0,0,0}$1,a$}%
}}}}
\put(3976,-1696){\makebox(0,0)[lb]{\smash{{\SetFigFont{8}{9.6}{\rmdefault}{\mddefault}{\updefault}{\color[rgb]{0,0,0}$1,c$}%
}}}}
\put(5326,-2221){\makebox(0,0)[lb]{\smash{{\SetFigFont{8}{9.6}{\rmdefault}{\mddefault}{\updefault}{\color[rgb]{0,0,0}$4,c$}%
}}}}
\end{picture}%

%% file: rec3_detail_alt.pstex_t
\begin{picture}(0,0)%
\includegraphics{rec3_detail_alt.pstex}%
\end{picture}%
\setlength{\unitlength}{2368sp}%
\begingroup\makeatletter\ifx\SetFigFont\undefined%
\gdef\SetFigFont#1#2#3#4#5{%
  \reset@font\fontsize{#1}{#2pt}%
  \fontfamily{#3}\fontseries{#4}\fontshape{#5}%
  \selectfont}%
\fi\endgroup%
\begin{picture}(4812,3174)(289,-2548)
\put(3901,-1936){\makebox(0,0)[lb]{\smash{{\SetFigFont{7}{8.4}{\rmdefault}{\mddefault}{\updefault}{\color[rgb]{0,0,0}$C_1$}%
}}}}
\put(5101,-1521){\makebox(0,0)[lb]{\smash{{\SetFigFont{9}{10.8}{\rmdefault}{\mddefault}{\updefault}{\color[rgb]{0,0,0}$\pi_3=\{C_1, C_2\}$}%
}}}}
\put(5101,-593){\makebox(0,0)[lb]{\smash{{\SetFigFont{9}{10.8}{\rmdefault}{\mddefault}{\updefault}{\color[rgb]{0,0,0}$\pi_1=\{A_1, A_2, A_3, A_4\}$}%
}}}}
\put(3376,239){\makebox(0,0)[lb]{\smash{{\SetFigFont{7}{8.4}{\rmdefault}{\mddefault}{\updefault}{\color[rgb]{0,0,0}$A_4$}%
}}}}
\put(2776,239){\makebox(0,0)[lb]{\smash{{\SetFigFont{7}{8.4}{\rmdefault}{\mddefault}{\updefault}{\color[rgb]{0,0,0}$A_3$}%
}}}}
\put(2176,239){\makebox(0,0)[lb]{\smash{{\SetFigFont{7}{8.4}{\rmdefault}{\mddefault}{\updefault}{\color[rgb]{0,0,0}$A_2$}%
}}}}
\put(1576,239){\makebox(0,0)[lb]{\smash{{\SetFigFont{7}{8.4}{\rmdefault}{\mddefault}{\updefault}{\color[rgb]{0,0,0}$A_1$}%
}}}}
\put(4126,-1186){\makebox(0,0)[lb]{\smash{{\SetFigFont{7}{8.4}{\rmdefault}{\mddefault}{\updefault}{\color[rgb]{0,0,0}$C_2$}%
}}}}
\put(1501,-2236){\makebox(0,0)[lb]{\smash{{\SetFigFont{7}{8.4}{\rmdefault}{\mddefault}{\updefault}{\color[rgb]{0,0,0}$B_3$}%
}}}}
\put(1126,-1711){\makebox(0,0)[lb]{\smash{{\SetFigFont{7}{8.4}{\rmdefault}{\mddefault}{\updefault}{\color[rgb]{0,0,0}$B_2$}%
}}}}
\put(751,-1111){\makebox(0,0)[lb]{\smash{{\SetFigFont{7}{8.4}{\rmdefault}{\mddefault}{\updefault}{\color[rgb]{0,0,0}$B_1$}%
}}}}
\put(5101,-1057){\makebox(0,0)[lb]{\smash{{\SetFigFont{9}{10.8}{\rmdefault}{\mddefault}{\updefault}{\color[rgb]{0,0,0}$\pi_2=\{B_1, B_2, B_3\}$}%
}}}}
\end{picture}%

%% file: arith_12.pstex_t
\begin{picture}(0,0)%
\includegraphics{arith_12.pstex}%
\end{picture}%
\setlength{\unitlength}{3947sp}%
\begingroup\makeatletter\ifx\SetFigFont\undefined%
\gdef\SetFigFont#1#2#3#4#5{%
  \reset@font\fontsize{#1}{#2pt}%
  \fontfamily{#3}\fontseries{#4}\fontshape{#5}%
  \selectfont}%
\fi\endgroup%
\begin{picture}(9038,1824)(139,-1123)
\put(8701,-736){\makebox(0,0)[lb]{\smash{{\SetFigFont{12}{14.4}{\rmdefault}{\mddefault}{\updefault}{\color[rgb]{0,0,0}$F_2$}%
}}}}
\put(451,-736){\makebox(0,0)[lb]{\smash{{\SetFigFont{12}{14.4}{\rmdefault}{\mddefault}{\updefault}{\color[rgb]{0,0,0}$F_1$}%
}}}}
\put(3901,-736){\makebox(0,0)[lb]{\smash{{\SetFigFont{12}{14.4}{\rmdefault}{\mddefault}{\updefault}{\color[rgb]{0,0,0}$F_2$}%
}}}}
\put(5251,-736){\makebox(0,0)[lb]{\smash{{\SetFigFont{12}{14.4}{\rmdefault}{\mddefault}{\updefault}{\color[rgb]{0,0,0}$F_1$}%
}}}}
\put(6586,-435){\makebox(0,0)[lb]{\smash{{\SetFigFont{8}{9.6}{\rmdefault}{\mddefault}{\updefault}{\color[rgb]{0,0,0}$X$}%
}}}}
\put(6961,376){\makebox(0,0)[lb]{\smash{{\SetFigFont{8}{9.6}{\rmdefault}{\mddefault}{\updefault}{\color[rgb]{0,0,0}$X$}%
}}}}
\put(7036,-885){\makebox(0,0)[lb]{\smash{{\SetFigFont{8}{9.6}{\rmdefault}{\mddefault}{\updefault}{\color[rgb]{0,0,0}$X$}%
}}}}
\put(6586,-135){\makebox(0,0)[lb]{\smash{{\SetFigFont{8}{9.6}{\rmdefault}{\mddefault}{\updefault}{\color[rgb]{0,0,0}$X$}%
}}}}
\put(2776,-865){\makebox(0,0)[lb]{\smash{{\SetFigFont{10}{12.0}{\rmdefault}{\mddefault}{\updefault}{\color[rgb]{0,0,0}$X$}%
}}}}
\put(2776,335){\makebox(0,0)[lb]{\smash{{\SetFigFont{10}{12.0}{\rmdefault}{\mddefault}{\updefault}{\color[rgb]{0,0,0}$X$}%
}}}}
\put(2776,-265){\makebox(0,0)[lb]{\smash{{\SetFigFont{10}{12.0}{\rmdefault}{\mddefault}{\updefault}{\color[rgb]{0,0,0}$X$}%
}}}}
\put(1576,-865){\makebox(0,0)[lb]{\smash{{\SetFigFont{10}{12.0}{\rmdefault}{\mddefault}{\updefault}{\color[rgb]{0,0,0}$X$}%
}}}}
\put(1576,-265){\makebox(0,0)[lb]{\smash{{\SetFigFont{10}{12.0}{\rmdefault}{\mddefault}{\updefault}{\color[rgb]{0,0,0}$X$}%
}}}}
\put(1576,335){\makebox(0,0)[lb]{\smash{{\SetFigFont{10}{12.0}{\rmdefault}{\mddefault}{\updefault}{\color[rgb]{0,0,0}$X$}%
}}}}
\put(6736,226){\makebox(0,0)[lb]{\smash{{\SetFigFont{8}{9.6}{\rmdefault}{\mddefault}{\updefault}{\color[rgb]{0,0,0}$X$}%
}}}}
\put(7411, 90){\makebox(0,0)[lb]{\smash{{\SetFigFont{8}{9.6}{\rmdefault}{\mddefault}{\updefault}{\color[rgb]{0,0,0}$X$}%
}}}}
\put(6811,-660){\makebox(0,0)[lb]{\smash{{\SetFigFont{8}{9.6}{\rmdefault}{\mddefault}{\updefault}{\color[rgb]{0,0,0}$X$}%
}}}}
\put(7411,-585){\makebox(0,0)[lb]{\smash{{\SetFigFont{8}{9.6}{\rmdefault}{\mddefault}{\updefault}{\color[rgb]{0,0,0}$X$}%
}}}}
\put(7186,-285){\makebox(0,0)[lb]{\smash{{\SetFigFont{8}{9.6}{\rmdefault}{\mddefault}{\updefault}{\color[rgb]{0,0,0}$X$}%
}}}}
\end{picture}%

%% file: molecular4.pstex_t
\begin{picture}(0,0)%
\includegraphics{molecular4.pstex}%
\end{picture}%
\setlength{\unitlength}{3947sp}%
\begingroup\makeatletter\ifx\SetFigFont\undefined%
\gdef\SetFigFont#1#2#3#4#5{%
  \reset@font\fontsize{#1}{#2pt}%
  \fontfamily{#3}\fontseries{#4}\fontshape{#5}%
  \selectfont}%
\fi\endgroup%
\begin{picture}(6399,5349)(664,-5773)
\put(4876,-886){\rotatebox{324.0}{\makebox(0,0)[lb]{\smash{{\SetFigFont{12}{14.4}{\rmdefault}{\mddefault}{\updefault}{\color[rgb]{0,0,0}$A$-orbits}%
}}}}}
\put(5645,-2094){\rotatebox{288.0}{\makebox(0,0)[lb]{\smash{{\SetFigFont{12}{14.4}{\rmdefault}{\mddefault}{\updefault}{\color[rgb]{0,0,0}$B$-orbits}%
}}}}}
\put(1693,-2815){\rotatebox{72.0}{\makebox(0,0)[lb]{\smash{{\SetFigFont{12}{14.4}{\rmdefault}{\mddefault}{\updefault}{\color[rgb]{0,0,0}$B$-orbits}%
}}}}}
\put(2326,-4636){\rotatebox{324.0}{\makebox(0,0)[lb]{\smash{{\SetFigFont{12}{14.4}{\rmdefault}{\mddefault}{\updefault}{\color[rgb]{0,0,0}$B$-orbits}%
}}}}}
\put(4876,-5011){\rotatebox{36.0}{\makebox(0,0)[lb]{\smash{{\SetFigFont{12}{14.4}{\rmdefault}{\mddefault}{\updefault}{\color[rgb]{0,0,0}$B$-orbits}%
}}}}}
\put(3411,-887){\makebox(0,0)[lb]{\smash{{\SetFigFont{12}{14.4}{\rmdefault}{\mddefault}{\updefault}{\color[rgb]{0,0,0}$B$-orbits}%
}}}}
\put(3451,-5461){\makebox(0,0)[lb]{\smash{{\SetFigFont{12}{14.4}{\rmdefault}{\mddefault}{\updefault}{\color[rgb]{0,0,0}$A$-orbits}%
}}}}
\put(6001,-4186){\rotatebox{72.0}{\makebox(0,0)[lb]{\smash{{\SetFigFont{12}{14.4}{\rmdefault}{\mddefault}{\updefault}{\color[rgb]{0,0,0}$A$-orbits}%
}}}}}
\put(1426,-3511){\rotatebox{288.0}{\makebox(0,0)[lb]{\smash{{\SetFigFont{12}{14.4}{\rmdefault}{\mddefault}{\updefault}{\color[rgb]{0,0,0}$A$-orbits}%
}}}}}
\put(2026,-1486){\rotatebox{36.0}{\makebox(0,0)[lb]{\smash{{\SetFigFont{12}{14.4}{\rmdefault}{\mddefault}{\updefault}{\color[rgb]{0,0,0}$A$-orbits}%
}}}}}
\end{picture}%

%% file: pi_sigma_beta.pstex_t
\begin{picture}(0,0)%
\includegraphics{pi_sigma_beta.pstex}%
\end{picture}%
\setlength{\unitlength}{3947sp}%
\begingroup\makeatletter\ifx\SetFigFont\undefined%
\gdef\SetFigFont#1#2#3#4#5{%
  \reset@font\fontsize{#1}{#2pt}%
  \fontfamily{#3}\fontseries{#4}\fontshape{#5}%
  \selectfont}%
\fi\endgroup%
\begin{picture}(2700,1956)(1651,-1999)
\put(2551,-1936){\makebox(0,0)[lb]{\smash{{\SetFigFont{14}{16.8}{\rmdefault}{\mddefault}{\updefault}{\color[rgb]{0,0,0}$\pi_{\sigma(2)}$}%
}}}}
\put(2701,-1036){\makebox(0,0)[lb]{\smash{{\SetFigFont{10}{12.0}{\rmdefault}{\mddefault}{\updefault}{\color[rgb]{0,0,0}$\beta_2$}%
}}}}
\put(4351,-1036){\makebox(0,0)[lb]{\smash{{\SetFigFont{10}{12.0}{\rmdefault}{\mddefault}{\updefault}{\color[rgb]{0,0,0}$\beta_k$}%
}}}}
\put(1801,-1036){\makebox(0,0)[lb]{\smash{{\SetFigFont{10}{12.0}{\rmdefault}{\mddefault}{\updefault}{\color[rgb]{0,0,0}$\beta_1$}%
}}}}
\put(1651,-211){\makebox(0,0)[lb]{\smash{{\SetFigFont{14}{16.8}{\rmdefault}{\mddefault}{\updefault}{\color[rgb]{0,0,0}$\pi_1$}%
}}}}
\put(2551,-211){\makebox(0,0)[lb]{\smash{{\SetFigFont{14}{16.8}{\rmdefault}{\mddefault}{\updefault}{\color[rgb]{0,0,0}$\pi_2$}%
}}}}
\put(4201,-211){\makebox(0,0)[lb]{\smash{{\SetFigFont{14}{16.8}{\rmdefault}{\mddefault}{\updefault}{\color[rgb]{0,0,0}$\pi_k$}%
}}}}
\put(4201,-1936){\makebox(0,0)[lb]{\smash{{\SetFigFont{14}{16.8}{\rmdefault}{\mddefault}{\updefault}{\color[rgb]{0,0,0}$\pi_{\sigma(k)}$}%
}}}}
\put(1651,-1936){\makebox(0,0)[lb]{\smash{{\SetFigFont{14}{16.8}{\rmdefault}{\mddefault}{\updefault}{\color[rgb]{0,0,0}$\pi_{\sigma(1)}$}%
}}}}
\end{picture}%

%% file: molecular_prime.pstex_t
\begin{picture}(0,0)%
\includegraphics{molecular_prime.pstex}%
\end{picture}%
\setlength{\unitlength}{3947sp}%
\begingroup\makeatletter\ifx\SetFigFont\undefined%
\gdef\SetFigFont#1#2#3#4#5{%
  \reset@font\fontsize{#1}{#2pt}%
  \fontfamily{#3}\fontseries{#4}\fontshape{#5}%
  \selectfont}%
\fi\endgroup%
\begin{picture}(3324,1224)(289,-523)
\end{picture}%

%% file: non-prime.pstex_t
\begin{picture}(0,0)%
\includegraphics{non-prime.pstex}%
\end{picture}%
\setlength{\unitlength}{3947sp}%
\begingroup\makeatletter\ifx\SetFigFont\undefined%
\gdef\SetFigFont#1#2#3#4#5{%
  \reset@font\fontsize{#1}{#2pt}%
  \fontfamily{#3}\fontseries{#4}\fontshape{#5}%
  \selectfont}%
\fi\endgroup%
\begin{picture}(4374,624)(889,-673)
\end{picture}%